%% file: new_version.tex
\documentclass[11pt,a4paper]{article}

\input{macros}

\title{Sum-of-Squares Hierarchy for the Gromov-Wasserstein Problem}

\author{Hoang Anh Tran, Binh Tuan Nguyen, Yong Sheng Soh}

\begin{document}

\maketitle

\begin{abstract}
The Gromov-Wasserstein (GW) problem is a variant of the classical optimal transport problem that allows one to compute meaningful transportation plans between incomparable spaces.  At an intuitive level, it seeks plans that minimize the discrepancy between metric evaluations of pairs of points.  The GW problem is typically cast as an instance of a non-convex quadratic program that is, unfortunately, intractable to solve.  In this paper, we describe tractable semidefinite relaxations of the GW problem based on the Sum-of-Squares (SOS) hierarchy.  We describe how the Putinar-type and the Schmüdgen-type moment hierarchies can be simplified using marginal constraints, and we prove convergence rates for these hierarchies towards computing global optimal solutions to the GW problem.  The proposed SOS hierarchies naturally induce a distance measure analogous to the distortion metrics, and we show that these are genuine distances in that they satisfy the triangle inequality.  In particular, the proposed SOS hierarchies provide computationally tractable proxies of the GW distance and the associated distortion distances (over metric measure spaces) that are otherwise intractable to compute.
\end{abstract}

\section{Introduction}

A metric measure space is a triplet $\mathbb{X}:=(\CX,d_{\CX},\mu_{\CX})$ whereby $(\CX,d_{\CX})$ specifies a metric space and where $\mu_{\CX}$ specifies a probability distribution over $\CX$.  Given a pair of metric measure spaces $\mathbb{X} = (\CX, d_{\CX}, \mu_{\CX})$ and $\mathbb{Y} = (\CY, d_{\CY}, \mu_{\CY})$, the {\em Gromov-Wasserstein} (GW) distance is defined as the solution to the following optimization instance
\begin{equation} \label{eq:gw}
\mathrm{GW}(\mathbb{X},\mathbb{Y}) := \underset{\pi \in \Pi (\mu_{\CX},\mu_{\CY})}{\inf} \int_{\CX \times \CY} \int_{\CX \times \CY} c(d_{\CX} (x_0,x_1), d_{\CY} (y_0,y_1)) ~ d \pi (x_0,y_0) d \pi (x_1,y_1)  .
\end{equation}
Here, $\Pi (\mu_{\CX},\mu_{\CY})$ denotes the set of all {\em couplings} of $\mu_{\CX}$ and $\mu_{\CY}$; that is, these are probability measures over $\CX \times \CY$ that satisfy the marginal constraints 
\begin{equation} \label{eq:marginal}
\int_{\CY} d \pi (x,y) = d \mu_{\CX}, \quad \int_{\CX} d \pi (x,y) = d \mu_{\CY}.
\end{equation}
Note: Because $\pi$ is a probability measure, we necessarily have $\pi \geq 0$.

The cost function $c(d_{\CX} (x_0,x_1), d_{\CY} (y_0,y_1))$ measures the difference between the metric spaces $(\CX,d_{\CX})$ and $(\CY,d_{\CY})$.  A prominent choice is to use
\begin{equation} \label{eq:lp_definition}
c(d_{\CX} (x_0,x_1), d_{\CY} (y_0,y_1)) := | d_{\CX} (x_0,x_1)^q - d_{\CY} (y_0,y_1)^q |^p,
\end{equation}
where $p,q \in [1, +\infty)$.  In this case, the resulting distance is specified as 
$$
\mathrm{GW}(\mathbb{X},\mathbb{Y})^{1/p}.
$$

The Gromov-Wasserstein distance in the form specified via \eqref{eq:gw} and \eqref{eq:lp_definition} was introduced early on in \cite{Mem:07}, and subsequently studied in greater detail in \cite{Mem:11,Sturm:23}.  The GW distance is closely related to the Gromov-Hausdorff distance \cite{Gromov:99}, as well as a GW-type distance in \cite{Sturm:06} that is based on embedding $\XX$ and $\YY$ into a common metric space rather than measuring differences in metric evaluations via \eqref{eq:lp_definition}.  The specific distance $\mathrm{GW}(\mathbb{X},\mathbb{Y})^{1/p}$ is also referred to as the $L^{p,q}$ {\em distortion distance} in \cite{Sturm:23}; in particular, the work in \cite{Sturm:23} introduces a series of geometric results about this distortion distance, such as geodesics and gradient flows over the {\em space of spaces}. The Gromov-Wasserstein distance has applications in several areas, such as graph matching \cite{Mem:11,xu2019gromov,vayer2020fused}, computational graphics \cite{solomon2016entropic,peyre2016gromov}, natural language processing \cite{alvarez2018gromov}, computational biology \cite{demetci2020gromov,klein2024genot}, brain imaging analysis \cite{thual2022aligning}, and generative modeling in machine learning \cite{bunne2019learning}.

{\bf Optimal transport distances over incomparable spaces.}  The Gromov-Wasserstein problem builds on the ideas of classical optimal transport.  Concretely, suppose one is given distributions $\mu_{\CX}$ and $\mu_{Y}$ over $\CX$ and $\CY$.  The (Kantorovich formulation of the) classical optimal transportation problem seeks an optimal coupling $\pi$, subject to the marginal constraints \eqref{eq:marginal}, so as to minimize the cost of moving the probability measure from $\mu_{\CX}$ to $\mu_{\CY}$:
\begin{equation} \label{eq:kantorovich_ot}
\inf_{\pi \in \Pi(\mu_X,\mu_Y)} ~~ \int_{\CX \times \CY} c(x,y) d \pi (x,y).
\end{equation}
Here, $c(x,y)$ denotes the cost of moving mass from location $x\in \CX$ to $y \in \CY$.

The formulation \eqref{eq:kantorovich_ot} pre-supposes that one has access to a cost function $c(x,y)$.  In certain applications where the spaces $\CX$ and $\CY$ are not naturally comparable, there is no natural notion of a cost $c(x,y)$ that one can define in order to obtain meaningful transportation maps.  This is frequently the case in applications such as shape matching across different dimensions and graph matching \cite{Mem:11}.  To this end, the Gromov-Wasserstein formulation \eqref{eq:gw} can be viewed as an alternative formulation that seeks to circumvent such a shortcoming by only requiring that one has knowledge of the metrics endowed on $\CX$ and $\CY$ respectively \cite{Mem:11,Sturm:06}.  At an intuitive level, the formulation in \eqref{eq:gw} can be viewed one that seeks a coupling $\pi$ that maximizes the {\em agreement} between the metrics endowed on $\CX$ and $\CY$ respectively.

\textbf{Computing the Gromov-Wasserstein distance.}  The GW problem as specified in \eqref{eq:gw} is defined over general probability distributions, including continuous ones.  In practical implementations, one typically discretizes the probability distributions $\muX$ and $\muY$ in order to compute an approximation of \eqref{eq:gw}.  We make such an assumption in this paper -- specifically, we assume that $\muX$ and $\muY$ are probability distributions that have {\em finite} support, say, on $\{x_1,\dots,x_m\} \subset \CX$ and $\{y_1,\dots,y_n\}$ respectively.  Under such an assumption, the GW problem \eqref{eq:gw} simplifies to the following non-convex {\em quadratic programming} instance:
\begin{equation} \label{eq:gw_discrete} \tag{DGW}
\underset{\pi}{\mathrm{min}} ~ \sum_{x_0,x_1,y_0,y_1} c(d_{\CX} (x_0,x_1), d_{\CY} (y_0,y_1)) ~ \pi_{x_0,y_0} \pi_{x_1,y_1} \quad \mathrm{s.t.} \quad \pi 1 = \mu_{\CX}, \pi^T 1 = \mu_{\CY}, \pi \geq 0.
\end{equation}

Unfortunately, this optimization instance is not known to be tractable to solve.  In particular, recent work in \cite{kravtsova2024np} shows that the general optimization formulation specified in \eqref{eq:gw_discrete} is in fact NP-Hard.  The result is perhaps not entirely surprising given the GW problem is known to be closely related to graph matching tasks and the Quadratic Assignment Problem (QAP), both of which are known to be NP-Hard \cite{bunne2019learning,conte2004thirty}.  The connection between these problems is as follows:  The graph matching task and the QAP may be stated as optimization instances whereby the objective is equal to that of \eqref{eq:gw_discrete} (the cost function $c$ needs not arise as a distortion measure between two metric spaces), but where the optimization variable $\pi$ is required to be a {\em permutation} matrix.  Recall from the Birkhoff-von Neumann theorem that the convex hull of permutation matrices is the set of doubly stochastic matrices.  In particular, the feasible region in \eqref{eq:gw_discrete} specialized to the case where $m=n$ specifies precisely the set of doubly stochastic matrices.  In some applications where graph matching tasks arise, one may sometimes view these graphs as arising via samples drawn from a continuous distribution -- in a sense, these graphs approximate the underlying continuous distribution.  

\subsection{Our contributions}

The objective of this paper is to provide a principled way of computing globally optimal transportation maps that minimize \eqref{eq:gw_discrete}.  We do so by applying the {\em sum-of-squares} (SOS) hierarchy (see e.g., \cite{lasserre2001global,lasserre2011new,BleParTho:12}), which is a principled approach for solving generic polynomial optimization instances via semidefinite relaxations of increasing size and approximation power to the original polynomial optimization instance.

Our first contribution is to state a reasonable implementation of the SOS hierarchy for the GW problem \eqref{eq:gw_discrete}.  This step is not as entirely trivial as those who are familiar with the SOS literature might think.  While the GW instance as it is stated in \eqref{eq:gw_discrete} specifies a polynomial optimization instance, it is not clear if this is the ideal formulation from the perspective of obtaining strong SOS-based relaxations.  As we discuss later in Section \ref{sec:moment-sos}, one has a choice between modelling non-negativity constraints on the variables $\pi_{ij}$ via linear inequalities, or by modelling these as squares $\pi_{ij} = \tilde{\pi}_{ij}^2$.  The latter choice seems unnecessarily complicated but actually leads to a description that is {\em Archimedean}, which is an important ingredient proving convergence of the hierarchy to the optimal value of \eqref{eq:gw_discrete}.  Our second contribution is to establish precisely such a convergence rate.

More importantly, these semidefinite relaxations are intended to be used as polynomial time-computable proxies or alternatives to the GW problem \eqref{eq:gw_discrete}, as a sub-routine to more complex computational tasks such as computing geodesics.  Our third contribution is to describe an alternative distance measure analogous to the $L^{p,q}$ distortion measure that is specified as the solution of a semidefinite program, and hence computable in polynomial time.  In particular, we show that these distance measures also satisfy the triangle inequality, and hence could in principle be used to compute geodesics with respect to the GW distance.  In fact, in prior work done in \cite{chen2023semidefinite}, the authors have noticed that semidefinite relaxations do indeed solve the discrete GW instance \eqref{eq:gw_discrete} to {\em global optimality} in many generic instances.  This raises the fascinating possibility of computing geodesics with respect to the SOS-analogs of the $L^{p,q}$ metrics in polynomial time that provably coincides with geodesics of the GW distance.

{\bf Prior work.}  We discuss prior work and our contributions in relation to these.  

First, the authors in \cite{villar2016polynomial,chen2023semidefinite} propose semidefinite relaxations to the (discretized)  GW problem \eqref{eq:gw_discrete}.  Our work generalize these ideas; in particular, the relaxation in \cite{chen2023semidefinite} can be viewed as the first level of the SOS hierarchy in our description.  On this note, we wish to point out that there is a substantive body of work that study semidefinite relaxations for the QAP \cite{zhao1998semidefinite,povh2009copositive,aflalo2015convex,kezurer2015tight,Ling:24,CS:24}.  In the same vein as the work in \cite{chen2023semidefinite}, these can be viewed as the first level of the SOS hierarchy whereas our goals are quite different.  (The QAP is also a different problem from GW.)

Second, the work in \cite{MulaSOS} investigates the use of SOS-hierarchy to a number of problems in optimal transport, including the Gromov-Wasserstein problem.  In this work, the authors suppose that the cost function -- this is represented by $c$ in \eqref{eq:kantorovich_ot} -- is a polynomial.  Based on this assumption, one applies the usual SOS hierarchy.  Our work makes a fundamentally different modelling choice -- namely -- that we view the parameters $\pi$ as the decision variables to be modelled as a polynomial while the cost function $c$ is to be modelled as coefficients of a polynomial.  There are several important consequences that result from these different modelling choices.  First, our work applies to general cost functions $c$ and not just polynomial ones.  Second, our framework naturally identifies the optimal transportation plan whereas the work in \cite{MulaSOS} are more suited for simply computing the actual objective value.  Third, and more generally,  our work can be viewed as the natural `continuous' analogs of semidefinite relaxations studied extensively in the context of QAPs.

Third, a piece of work that has some conceptual connections is that of \cite{ZGMS:24}.  In Section \ref{sec:moment-sos} we discuss primal and dual formulations of the SOS hierarchy for \eqref{eq:gw_discrete}.  The work of \cite{ZGMS:24} also introduces some notion of duality, though it is very different from ours.  In particular, the notion of duality in \cite{ZGMS:24} heavily relies on cost to take a specific polynomial form for which certain algebraic manipulations are possible.  Also, the notion of duality we develop is precisely Lagrangian duality whereas the duality in \cite{ZGMS:24} is quite different.

\section{Background}

The GW task \eqref{eq:gw_discrete} is a quadratic polynomial optimization instance over a $mn$-dimensional polytope.  Polynomial optimization problems (POPs) are intractable optimization problems in general; in particular, the optimization formulation we study in \eqref{eq:gw_discrete} is NP-hard, as we noted earlier.

To address the difficult problem of solving POP instances, the moment-sum of squares (SOS) hierarchy introduces a series of semidefinite relaxations that increasingly approximate the solution of the original POP instance \cite{lasserre2001global,thesisparrilo}.  Each optimization instance within the hierarchy is specified as a semidefinite program -- an optimization instance in which one minimizes a linear function over a constraint set specified as the intersection of the cone of positive semidefinite (PSD) matrices with an affine space -- and is solvable in polynomial time \cite{vandenberghe1996semidefinite,Ren:01}.

We begin by providing a brief background to the {\em moment}-SOS hierarchy for general POP instances.  First, given a multi-index $\alpha \in \NN^n$, we denote the associated monomial in $x$ by $x^{\alpha} := \prod_{i=1}^n x_i^{\alpha_i}$.  
Let $\CX \in \mathbb{R}^n$ be a semi-algebraic set specified as follows:
\begin{equation*}\label{def of semialg set}
\CX := \{x \in \mathbb{R}^n\,:\, g_j(x) \geq 0 \quad \forall j \in [m]\}.  
\end{equation*} 
Here, $g_j \in \mathbb{R}[x]$ are polynomials.  In what follows, we consider the POP instance in which we minimize a polynomial $f(x) := \sum_{\alpha \in \NN^n}f_{\alpha}x^{\alpha}$ over the constraint set $\CX$:
\begin{equation}\label{pop}
f_{\min} ~:=~ \min \{ f(x) : x \in \CX \}.
\end{equation}

Second, we let $\Sigma[x]$ denote the set of SOS polynomials
\begin{equation*}
\Sigma [x] ~:=~ \Big \{ \sum_{i=1}^k p_i^2 : p_i \in \RR[x],\ k \in \NN \Big \}.
\end{equation*}
In a similar fashion, we let $ \Sigma[x]_{2r} \subset \Sigma [x]$ denote the subset of SOS polynomials with degree at most $2r$.  The {\em quadratic module} $\CQ(\CX)$ and the {\em pre-ordering} $\CT(\CX)$ are subsets of $\mathbb{R}[x]$ defined as follows:
\begin{equation*}
\begin{aligned}
\CQ(\CX) ~&:=~ \Big\{ \sigma_0(x) +\sum_{i=1}^m\sigma_i(x)g_i(x) \in \RR[x] \,:\, \sigma_i \in \Sigma[x] \; \forall 0 \leq i \leq m\Big\},\\
\CT(\CX) ~&:=~ \Big\{\sum_{j=1}^N\sigma_{I_j}(x)g_{I_j}(x) \in \RR[x] \,:\, N \in \NN,\ I_j \subset [m],\ \sigma_{I_j} \in \Sigma[x]\ \Big\},
\end{aligned}
\end{equation*}
where for any subset $I \subset [m]$, we define $g_I = \prod_{i \in I}g_i$ with the convention $ g_{\emptyset}=1.$
Note that these specify convex cones.  Let $d_i := \lceil \deg g_i /2 \rceil$.  Then the corresponding degree-bounded analogs of the quadratic module and the pre-ordering are given by
\begin{equation*}
\begin{aligned}
\CQ(\CX)_{2r} ~&:=~ \Big\{\sigma_0(x) +\sum_{i=1}^m\sigma_i(x)g_i(x) \in \CQ(\CX) \,:\, \sigma_0 \in \Sigma[x]_{2r},\; \sigma_i \in \Sigma[x]_{2r-2d_i}\; \forall 0 \leq i \leq m  \Big\},\\
\CT(\CX)_{2r} ~&:=~ \Big\{\sum_{j=1}^N\sigma_{I_j}(x)g_{I_j}(x) \in \CT(\CX) \,:\, \deg \sigma_{I_j}g_{I_j} \leq 2r \ \forall j \in [N] \Big\}.
\end{aligned}
\end{equation*}

Third, given a sequence of reals $\iy = (\iy_{\alpha})_{\alpha \in \NN^n} \subset \RR$ indexed by monomials in $x$, the {\em moment matrix} $\M(\iy)$ is defined by
\begin{equation*}
    \M(\iy)(\alpha, \beta) =\iy_{\alpha +\beta}, \quad \text{ for all } \quad \alpha, \beta \in \NN^n.
\end{equation*}
Given an integer $r \in \NN$, the $r$-truncated moment matrix $\M_r(\iy)$ is the symmetric matrix of dimension $s(n,r) := \binom{n+r}{r}$ obtained by restricting the rows and columns $\M(\iy)$ to monomials of degree at most $r$:
\begin{equation*}
    \M(\iy)(\alpha, \beta) =\iy_{\alpha +\beta}, \quad \forall \alpha, \beta \in \NN^n_r,
\end{equation*}
where $\NN^n_r$ is the set of multi-indices with length at most $r$; i.e., $\NN^n_r := \{\alpha = (\alpha_1,\dots,\alpha_n) \in \NN^n :  |\alpha| := \sum_{i=1}^n \alpha_i \leq r\}$.
The moment matrices and the localizing matrices are used to describe conditions under which a sequence of reals $\iy$ can be realized as a sequence of moments corresponding to some Borel measure (see for instance the Riesz-Haviland theorem, e.g., \cite[Theorem 3.1]{lasserre2009moments}).  We extend these definitions to allow more general inputs.  Specifically, given a polynomial $g(x) = \sum_{\gamma \in \NN^n}g_{\gamma}x^{\gamma} \in \RR[x]$, the {\em localizing matrix} $\M(g\iy)$ associated with the sequence $\iy$ and the polynomial $g$ is defined by 
\begin{equation*}
    \M(g\iy)(\alpha,\beta)=\sum_{\gamma \in \NN^n}g_{\gamma}\iy_{\gamma+\alpha+\beta}, \quad \text{ for all } \quad \alpha, \beta \in \NN^n.
\end{equation*}
In a similar fashion, $\M_r(g\iy)$ denotes the $r$-truncated localizing matrix obtained by restricting the rows and columns of $\M(g\iy)$ to monomials of degree at most $r$. 

We are now in a position to state the relevant moment-SOS hierarchies.  Concretely, the $r$-th level of the {\em Schmüdgen-type} moment hierarchy is given by
\begin{equation*}
\begin{aligned}
\mlb(f,\mathcal{T}(\CX))_r ~:=~ \underset{\iy \in \RR^{s(n,2r)}}{\inf} \quad & \sum_{\alpha \in \NN_{2r}^n} f_{\alpha}\iy_{\alpha}\\
\text{subject to } \quad & \iy_0=1,\; \M_{r-d_I}(g_I\iy) \succeq 0 \, \forall \  \{ I : I \subset [m], \deg g_I \leq 2r \}
\end{aligned}.
\end{equation*}
Here, we define $d_I := \lceil \deg g_I /2 \rceil$.  Note that the (constant) polynomial $1$ has degree $0$ so the above constraints automatically include the constraint $\M_r(\iy) \succeq 0$ by definition.  A consequence of the {\em Schmüdgen Positivstellensatz} is that, in the event where $\CX$ is compact, the sequence of lower bounds $\{\mlb(f,\mathcal{T}(\CX))_r\}_{r}^{\infty}$ increases monotonically to $f_{\min}$ as $r \rightarrow \infty$ \cite{schmudgen2017moment}.  (Generally speaking, {\em Positivstellensätze} are results that provide conditions under which positivity of a polynomial over certain sets can be certified via SOS polynomials.)

In a similar fashion, the $r$-th level of the {\em Putinar-type} moment hierarchy is given by 
\begin{equation*}
\begin{aligned}
\mlb(f,\mathcal{Q}(\CX))_r ~:=~ \underset{\iy \in \RR^{s(n,2r)}}{\inf} \quad & \sum_{\alpha \in \NN_{2r}^n} f_{\alpha}\iy_{\alpha} \\
\text{subject to } \quad & \iy \in \RR^{s(n,2r)},\; \iy_0=1,\; \M_r(\iy) \succeq 0,\; \M_{r-d_i}(g_i\iy) \succeq 0 \ \forall  i \in [m]
\end{aligned}.
\end{equation*}
We say that the constraint set $\CX$ satisfies the {\em Archimedean condition} if there exists $R>0$ such that 
\begin{equation*}
    R - \|x\|^2 \in \CQ(\CX). 
\end{equation*}
A consequence of the {\em Putinar Positivstellensatz} is that, in the event where $\CX$ is Archimedean, the $\{\mlb(f,\mathcal{Q}(\CX))_r\}_{r}^{\infty}$ also monotonically increases to $f_{\min}$ as $r \rightarrow \infty$ \cite{putinar1993positive}.

These hierarchies have very different characteristics:  The Schmüdgen-type hierarchy requires us to specify constraints arising from the localizing matrices generated from {\em all possible products} of the polynomials that specify the constraint set $\CX$; in contrast, the Putinar-type hierarchy only requires products of a {\em single} term, and hence leads to SDPs with far more compact descriptions.  On the other hand, where providing conditions that guarantee convergence is concerned, the Putinar-type hierarchy typically requires the constraint set to be Archimedean whereas Schmüdgen-type hierarchy only requires compactness, which is weaker.


{\bf Dual formulations and interpretations.}  We describe the (Lagrangian) dual formulation of the moment-SOS hierarchies -- these are known as the {\em SOS hierarchy}, and have the natural interpretation of linear optimization instances over the cone of SOS polynomials (with the appropriate degree bounds).  The Schmüdgen-type SOS hierarchy is given by

\begin{equation*}
\lb(f,\CT(\CX))_r = \sup\{c \in \RR:\ f(x) -c \in \CT(\CX)_{2r}\},
\end{equation*}
while the Putinar-type SOS hierarchy is given by 
\begin{equation*}
\lb(f,\CQ(\CX))_r = \sup\{c \in \RR:\ f(x) -c \in \CQ(\CX)_{2r}\}.
\end{equation*}

\section{Moment-SOS hierarchy for Discrete GW} \label{sec:moment-sos}

In this section we describe the moment-SOS hierarchy for the discrete GW problem \eqref{eq:gw_discrete}.  First we write $\Pi(\mu_{\CX},\mu_{\CY})$ as follows
\begin{equation*}
\begin{aligned}
\Pi(\mu_{\CX},\mu_{\CY}) = \bigl\{ \pi \in \RR^{m \times n}~:~ & e_{ij}(\pi) := \pi_{ij} \geq 0 ~ \forall i \in [m],\; j \in [n], \\
& r_i(\pi) := \textstyle \sum_{j \in [n]} \pi_{ij} - \muxi =0 ~ \forall i \in [m], \\
& c_j(\pi):= \textstyle \sum_{i \in [m]} \pi_{ij} - \muyj =0 ~ \forall j \in [n] \bigr\}.
\end{aligned}
\end{equation*}
Here, $\Pi(\mu_{\CX},\mu_{\CY})$ is a semialgebraic set defined by $mn$ linear inequalities of the form $e_{ij}(\pi) \geq 0$, and $m+n$ affine equalities of the form $r_i(\pi) = 0$ and $c_j(\pi)=0$.  The problem \eqref{eq:gw_discrete} can then be rewritten in the following form
\begin{equation}\tag{DGW+}\label{eq:gw_discrete+}
    \min_{\pi \in \Pi(\mu_{\CX},\mu_{\CY})}\sum_{i \in [m], j\in [n]}L_{ij,kl}\cdot \pi_{ij}\pi_{kl},
\end{equation}
where $L_{ij,kl} = c(d_{\CX}(x_i,x_k),d_{\CY}(y_j,y_l))$.
It is not clear that the set $\Pi(\mu_{\CX},\mu_{\CY})$ as it is specified satisfies the Archimedean property, and this suggests that the Schmüdgen-type moment-SOS hierarchy is the relevant procedure to apply if one wishes to develop a hierarchy that provably converges to the GW problem. 

As we noted, the main downside of the Schmüdgen-type hierarchy is the substantial number of constraints.  A direct application of the Schmüdgen-type hierarchy at the $r$-th level leads to the following semidefinite relaxation comprising $\sum_{k=0}^{2r}\binom{mn+m+n}{k}$ constraints
\begin{equation} \label{eq:s-type-constraints}
\M_{r-\lceil \deg g /2\rceil}(g \iy) * 0,\; g = \prod_{i=1}^k g_i,\; \{g_1,\dots,g_k\} \subset \{e_{ij},r_{i},c_j\},\; * = 
\begin{cases}
\succeq & \text{if } \{g_1,\dots,g_k\} \subset \{e_{ij}\}\\
= & \text{otherwise }
\end{cases}.
\end{equation}
This is a substantial number of constraints, and especially so for higher levels of the moment-SOS hierarchy!  In what follows, we derive an alternative scheme based on the Putinar-type hierarchy that still provably converges.

{\bf Eliminating redundant constraints.}  First, we make a simple observation that allows us to eliminate a number of redundant constraints in \eqref{eq:s-type-constraints}.  Consider the sequence $\iy = (\iy_{\alpha})_{\alpha \in \NN^{mn}_{2r}} \subset \RR^{s(mn,2r)}$ whereby $\NN^{mn}_{2r}$ denotes the set of multi-index $\alpha \in \NN^{mn}$ of length $|\alpha| = \sum_{i=1}^{mn}\alpha_i \leq 2r$. 
To simplify notation, we denote the component of $\iy$ indexed by the monomial $\pi_{ij}$ with $\iy_{ij}$, and the component of $\iy$ indexed by the monomial $\pi_{ij}\pi_{kl}$ with $\iy_{ij,kl}$.  In what follows, we let $[ m \times n]:= \{(i,j):\; i \in [m],\; j \in [n]\}$ denote the set indices within a $m \times n-$size matrix $\RR^{m \times n}$. Given a subset of indices $I \subset [m \times n]$, we denote the multi-index $\gamma_I \in \NN^{mn}$ whereby the ${ij}$-th component is $1$ if $ij \in I$, and is $0$ otherwise.  Subsequently, we let $e_I$ denote the monomial $\pi^{\gamma_I}:=\prod_{ij \in I} \pi_{ij}$, and we let $d_I := \lceil \deg e_I /2 \rceil $ denote the associated degree.  Consider the following version of the Schmüdgen-type moment hierarchy whereby the $r$-th level of the hierarchy is given as follows
\begin{equation} \label{hierarchy Schmudgen discrete GW} 
\begin{aligned}
\min \qquad & \sum_{i,k \in [m],\; j,l \in [n]}L_{ij,kl}\cdot \iy_{ij,kl} \\
\mbox{subject to} \qquad & \iy \in \RR^{s(mn,2r)},\;\iy_0 =1, \;\\ 
& \M_{r - d_I}(e_I\iy) \succeq 0 \; \forall I \subset [m \times n],\; \deg e_I \leq 2r,\\
& \M_{r-1}(r_i(\pi) \iy) = 0, \;\M_{r-1}(r_{i_1}e_{i_2j}(\pi) \iy) =0\; \forall i_1, i_2 \in [m],\ j \in [n]\\
& \M_{r-1}(c_j(\pi) \iy) = 0, \;\M_{r-1}(c_{j_1}e_{ij_2}(\pi) \iy) =0\; \forall j_1,j_2 \in [n],\ i \in [m].
\end{aligned} \tag{S-DGW-r}
\end{equation}

We explain why the constraints in \eqref{hierarchy Schmudgen discrete GW} imply the constraints in \eqref{eq:s-type-constraints}.  Consider the following subset of conditions
\begin{equation}\tag{Mar}\label{marginal conditions}
\begin{aligned}
& \M_{r-1}(r_i(\pi) \iy) = 0, \;\M_{r-1}(r_{i_1}e_{i_2j}(\pi) \iy) =0\; \forall i_1, i_2 \in [m],\ j \in [n]\\
&\M_{r-1}(c_j(\pi) \iy) = 0, \;\M_{r-1}(c_{j_1}e_{ij_2}(\pi) \iy) =0\; \forall j_1,j_2 \in [n],\ i \in [m]
\end{aligned}
\end{equation}
The $(0,0)$-th coordinate of these matrices enforce the constraints $r_i(\pi)=0$ and $c_j(\pi)=0$.  These in turn imply that any PSD constraint associated to products $g$ of polynomials $e_{ij}$'s, $r_i$'s, and $c_j$'s that contain at least one copy of $r_i$ or $c_j$ in \eqref{eq:s-type-constraints} necessarily hold -- in fact, they hold with {\em equality}.  In particular, the constraints in \eqref{hierarchy Schmudgen discrete GW} imply those in \eqref{eq:s-type-constraints}.  

Unfortunately, even with the reductions due to the inclusion of \eqref{marginal conditions}, the resulting semidefinite relaxation \eqref{hierarchy Schmudgen discrete GW} still contains $\sum_{k=0}^{2r}\binom{mn}{k}$ PSD constraints, which is prohibitive.  This necessitates the ideas in the next section.

{\bf Reduced Putinar-type hierarchy.}  In the following, we describe an alternative parameterization of the GW problem whereby the constraint set satisfies the Archimedean condition.  Concretely, we perform the following mapping
\begin{equation*}
    \varphi: \RR^{mn} \to \RR^{mn}: \; (\pi_{ij}) \mapsto (\widetilde{\pi}_{ij}^2).
\end{equation*}
Recall that the variables $\pi_{ij}$ model probabilities and hence are required to be non-negative.  As such, it is valid to model these as squares.  It is clear that the pre-image of $\widetilde{\Pi}(\mu_{\CX},\mu_{\CY})$ is contained in the $mn$-dimensional unit ball $B^{mn}$; in particular, one has
%
\begin{equation*}
    \widetilde{\Pi}(\mu_{\CX},\mu_{\CY}) = \big \{(\widetilde{\pi}_{ij})_{m \times n} \in \RR^{mn}:\; \textstyle \sum_{j\in [n]}\widetilde{\pi}_{ij}^2 = \muxi,\; \textstyle \sum_{i\in [m]}\widetilde{\pi}_{ij}^2 = \muyj\; \forall i \in [m],\; j \in [n] \big\}.
\end{equation*}
This implies that 
\begin{equation*}
   1 -\left\|(\widetilde{\pi}_{ij})_{m \times n}\right\|^2 = \textstyle \sum_{i \in [m]}\big(\muxi-\textstyle \sum_{j\in [n]}\widetilde{\pi}_{ij}^2 \big) \in \CQ(\widetilde{\Pi}(\mu_{\CX},\mu_{\CY}));
\end{equation*}
that is, the set $\widetilde{\Pi}(\mu_{\CX},\mu_{\CY})$ satisfies the Archimedean condition (we pick $\sigma_i=1$ and $g_i = \muxi-\sum_{j}\widetilde{\pi}_{ij}^2$ in the description).  The mapping $\varphi$ induces a pullback on the problem \eqref{eq:gw_discrete}; namely, \eqref{eq:gw_discrete} is equivalent to the following POP
\begin{equation}\label{discrete GW with pull back}
\eqref{eq:gw_discrete} ~=~ \min \Big \{\widetilde{L}(\widetilde{\pi}):= \textstyle \sum_{i,k \in [m], j,l \in [n]}L_{ij,kl}\cdot\widetilde{\pi}_{ij}^2\widetilde{\pi}_{kl}^2 ~:~ \widetilde{\pi} \in \widetilde{\Pi}(\mu_{\CX},\mu_{\CY}) \Big \}.
\end{equation}
This is a biquadratic polynomial optimization instance whereby the constraint set $\widetilde{\Pi}(\mu_{\CX},\mu_{\CY})$ is semialgebraic and satisfies the Archimedean condition. (A biquadratic function is quartic polynomial with no odd powers.)  In particular, we are now able to apply the Putinar-type moment-SOS hierarchy to obtain a hierarchy that provably converges to \eqref{eq:gw_discrete}.

In what follows, we apply a similar simplification step analogous to \eqref{hierarchy Schmudgen discrete GW}.  More specifically, every sequence $(\widetilde{\iy}_{\alpha})_{\alpha \in \NN^{mn}_{2r}} \subset \RR^{s(mn,2r)}$ is indexed by multi-indices of length at most $2r$, each of which corresponds to a monomial of degree at most $2r$.  And as such, with a slight abuse of notation, we simply say that the sequence $(\widetilde{\iy}_{\alpha})_{\alpha \in \NN^{mn}_{2r}}$ is indexed by monomials of degree at most $2r$.   In what follows, we denote the component of $\widetilde{\iy}$ indexed by the monomials $\widetilde{\pi}_{ij}^2$ with $\widetilde{\iy}_{ij,ij}$ and the component of $\widetilde{\iy}$ indexed by the monomial $\widetilde{\pi}_{ij}^2\widetilde{\pi}_{kl}^2$ with $\widetilde{\iy}_{ij,ij,kl,kl}$. The proposed Putinar-type moment hierarchy is thus given by
\begin{equation} \label{hierarchy Putinar discrete GW}
\begin{aligned}
\min \qquad & \widetilde{L}(\widetilde{\pi}):=\sum_{i,k \in [m],\; j,l \in [n]}L_{ij,kl}\cdot \widetilde{\iy}_{ij,ij,kl,,kl}\\
\mbox{subject to} \qquad & \widetilde{\iy} \in \RR^{s(mn,2r)}, \;\widetilde{\iy}_0 =1,\;\M_r(\widetilde{\iy}) \succeq 0, \\ 
& \M_{r-1}((\muxi - \textstyle\sum_{j \in [n]}\widetilde{\pi}_{ij}^2)\widetilde{\iy}) = 0 \; \forall i \in [m],\\
&  \M_{r-1}((\muyj - \textstyle\sum_{i \in [m]}\widetilde{\pi}_{ij}^2)\widetilde{\iy}) = 0 \; \forall j \in [n].
\end{aligned}
\tag{P-DGW-r}
\end{equation}


Notice that there is only {\em one} semidefinite constraint and $m+n$ linear constraints, which is substantially fewer than the number of constraints of the Schmüdgen-type moment hierarchy of \eqref{eq:gw_discrete} and \eqref{hierarchy Schmudgen discrete GW}.  There is a cost --unfortunately, the mapping $\varphi$ doubles the degree of the relaxation.  Our next result formalizes the connection between these hierarchies.

\begin{theorem}\label{connection between 2 hierarchies}
Let $r$ be any positive integer.  One has $\eqref{hierarchy Putinar discrete GW of double degree} = \eqref{hierarchy Schmudgen discrete GW} \leq \eqref{eq:gw_discrete}$. 
\end{theorem}

\begin{proof}[Proof of Theorem \ref{connection between 2 hierarchies}]
The inequality $\eqref{hierarchy Schmudgen discrete GW} \leq \eqref{eq:gw_discrete}$ follows from the fact that each level of the hierarchy $\{ \eqref{hierarchy Schmudgen discrete GW}\}_{r}^{\infty}$ provides a lower bound to \eqref{eq:gw_discrete}.  It remains to prove the equality $\eqref{hierarchy Putinar discrete GW of double degree} = \eqref{hierarchy Schmudgen discrete GW}$.  

Recall that 
\begin{equation} \label{hierarchy Putinar discrete GW of double degree}
\begin{aligned}
\min \qquad & \sum_{i,k \in [m],\; j,l \in [n]}L_{ij,kl}\cdot \widetilde{\iy}_{ij,ij,kl,,kl}\\
\mbox{subject to} \qquad & \widetilde{\iy} \in \RR^{s(mn,4r)},\; \widetilde{\iy}_0 =1,\; \M_{2r}(\widetilde{\iy}) \succeq 0,\\ 
& \M_{2r-1}((\muxi - \textstyle\sum_{j \in [n]}\widetilde{\pi}_{ij}^2)\widetilde{\iy}) = 0 \; \forall i \in [m],\\
& \M_{2r-1}((\muyj - \textstyle\sum_{i \in [m]}\widetilde{\pi}_{ij}^2)\widetilde{\iy}) = 0 \; \forall j \in [n].
\end{aligned}\tag{P-DGW-2r}
\end{equation}

[$\eqref{hierarchy Schmudgen discrete GW} \leq \eqref{hierarchy Putinar discrete GW of double degree}$]:  Given $\widetilde{\iy} \in \RR^{s(mn,4r)}$, we define a projection $\Proj: \RR^{s(mn,4r)} \to \RR^{s(mn,2r)}$ specified by 
\begin{equation*}
     \Proj(\widetilde{\iy})_{\alpha} = \widetilde{\iy}_{2\alpha} \quad \forall \; \widetilde{\iy}=(\widetilde{\iy}_{\gamma})_{\gamma \in \NN^{mn}_{4r}} \in \RR^{s(mn,4r)},\;  \alpha \in \NN^{mn}_{2r}.
\end{equation*} 
Suppose that $\widetilde{\iy}$ is a feasible solution to \eqref{hierarchy Putinar discrete GW of double degree}.  

We begin by showing that the projection $\Proj(\widetilde{\iy})$ also specifies a feasible solution to \eqref{hierarchy Schmudgen discrete GW}.  First, we note that marginal conditions of \eqref{hierarchy Schmudgen discrete GW} are straightforward to verify and we omit these steps.  We proceed to show that $\M_{r - d_I}(e_I P(\widetilde{\iy})) \succeq 0$ for all $I \subset [m \times n]$, $\deg e_I \leq 2r$.

We recall a few facts:  Given an integer $r$, let $\RR[\pi]_{r}$ denote the space of polynomials with degree at most $r$.  Given a polynomial $p \in \RR[\pi]_{r}$, one can represent $p(\pi) = \sum_{\alpha \in \NN^{mn}_r}p_{\alpha}\pi^{\alpha}$ via its vector of coefficients $(p_{\alpha})_{\alpha \in \NN^{mn}_r}$.  The inner product associated with $\M_r(\iy) \succeq 0$, where $\iy \in \RR^{s(mn,r)}$, is given by $ \langle p,q \rangle_{\iy} = \sum_{\alpha, \beta \in \NN^{mn}_r}\iy_{\alpha+\beta}p_{\alpha}q_{\beta}$.  

Based on the above, showing that $\M_{r - d_I}(e_I P(\widetilde{\iy})) \succeq 0$ is equivalent to showing that 
\begin{equation*}
    \langle p,e_Ip \rangle_{\Proj(\widetilde{\iy})} = p^{\top}\M_{r- \lceil \deg e_I/2 \rceil}(e_I \Proj(\widetilde{\iy}))p \geq 0
\end{equation*}
for any polynomial $p$ such that $2\deg p + \deg e_I \leq 2r$.  Since the $\alpha$-component of the projection $\Proj(\widetilde{\iy})$ is the $2\alpha$-component of the sequence $\widetilde{\iy}$, while $e_I$ is a product of monomials in $\pi$, for any $p, q \in \RR[\pi]_r$, one can represent $\langle p,q \rangle_{\Proj(\widetilde{\iy})}$ in terms of $\langle \cdot,\cdot \rangle_{\widetilde{\iy}}$ as 
\begin{equation*}
    \langle p,q \rangle_{\Proj(\widetilde{\iy})}= \langle p \circ \varphi, q \circ \varphi \rangle_{\widetilde{\iy}}, \quad \varphi \;:\; \RR^{mn} \to \RR^{mn} \;:\; \pi \mapsto \widetilde{\pi}^2.
\end{equation*}
Subsequently, let $p$ be any polynomial such that $2\deg p + \deg e_I \leq 2r$.  Then one has
\begin{equation*}
    \langle p,e_Ip \rangle_{\Proj(\widetilde{\iy})} = \langle p \circ \varphi, e_I \circ  \varphi \cdot p \circ \varphi \rangle_{\widetilde{\iy}}.
\end{equation*}
Since $e_I$ is a monomial, one has $e_I \circ \varphi = e_I^2 $, from which we obtain 
\begin{equation*}
\langle p,e_Ip \rangle_{\Proj(\widetilde{\iy})} = \langle p \circ \varphi, e_I^2 \cdot p \circ \varphi \rangle_{\widetilde{\iy}} = \langle e_I \cdot p \circ \varphi, e_I\cdot p \circ \varphi \rangle_{\widetilde{\iy}} = (e_I\cdot p \circ \varphi)^{\top}\M_{2r}(\widetilde{\iy})(e_I\cdot p \circ \varphi) \geq 0.
\end{equation*}
The last inequality is based on the fact that the matrix $\M_{2r}(\widetilde{\iy})$ is positive semidefinite. 

This concludes the proof that the projection of any feasible solution $\widetilde{\iy}$ of \eqref{hierarchy Putinar discrete GW of double degree} via $\Proj$ is a feasible solution of \eqref{hierarchy Schmudgen discrete GW}.  Finally, we check that 
\begin{equation*}
\textstyle \sum_{i,k \in [m],\; j,l \in [n]}L_{ij,kl}\cdot \widetilde{\iy}_{ij,ij,kl,,kl}= \textstyle \sum_{i,k \in [m],\; j,l \in [n]}L_{ij,kl}\cdot \Proj(\widetilde{\iy})_{ij,kl}.
\end{equation*}
By minimizing over $\widetilde{\iy}$, we obtain the conclusion $\eqref{hierarchy Schmudgen discrete GW} \leq \eqref{hierarchy Putinar discrete GW of double degree}$.

[$\eqref{hierarchy Putinar discrete GW of double degree} \leq \eqref{hierarchy Schmudgen discrete GW}$]:  Next we prove the reverse direction.  To do so, we show that every feasible solution of \eqref{hierarchy Schmudgen discrete GW} $\iy$ can be extended to be a feasible solution $\widetilde{\iy}$ of \eqref{hierarchy Putinar discrete GW of double degree}.  Concretely, consider the mapping $\mathbf{Q}: \RR^{s(mn,2r)} \to \RR^{s(mn,4r)}$ defined by 
\begin{equation*}
\widetilde{\iy}:=\mathbf{Q}(\iy)_{2\alpha} = \iy_{\alpha} \; \forall \alpha \in \NN^{mn}_{2r}, \quad \mathbf{Q}(\iy)_{\beta} = 0 ~ \text{otherwise}.
\end{equation*}
It is straightforward from the definition that for all feasible solution $\iy$ of \eqref{hierarchy Schmudgen discrete GW}, $\widetilde{\iy}$ defined via $\varphi$ satisfies the linear constraints associated to the marginal conditions
\begin{equation*}
\M_{r-1}((\muxi - \textstyle \sum_{j \in [n]}\widetilde{\pi}_{ij}^2)\widetilde{\iy}) = 0 \; \forall i \in [m], \quad \M_{r-1}((\muyj - \textstyle \sum_{i \in [m]}\widetilde{\pi}_{ij}^2)\widetilde{\iy}) = 0 \; \forall j \in [n].
\end{equation*}

It remains to prove that $\widetilde{\iy}$ satisfies the semidefinite constraint $\M_{2r}(\widetilde{\iy}) \succeq 0$.  In particular, one needs to show that, for any polynomial $p(\widetilde{\pi}) \in \RR[\widetilde{\pi}]_{2r}$, one also has $\langle p, p \rangle_{\widetilde{\iy}} \geq 0$.  We expand the product $\langle p, p \rangle_{\widetilde{\iy}}$ to obtain
\begin{equation*}
    \langle p, p \rangle_{\widetilde{\iy}} = \textstyle \sum_{\alpha, \beta \in \NN^{mn}_{2r}}\widetilde{\iy}_{\alpha+\beta}p_{\alpha}p_{\beta}.
\end{equation*}
Notice that $\widetilde{\iy}_{\alpha + \beta}$ vanishes unless $\alpha+\beta = 2 \gamma$ for some $\gamma$.  Therefore we only take into account pairs of multi-indices of the form $\alpha = 2\alpha^\prime + \gamma_I$ and $\beta = 2\beta^\prime + \gamma_I$ for some $\alpha^\prime, \beta^\prime \in \NN^{mn}_r$, and for some $I \subset [m \times n]$, where $\gamma_I$ represents the odd part of the multi-index $\alpha$ and $\beta$. As such, we can partition the above total sum as follows:
\begin{equation*}
\begin{aligned}
& \textstyle \sum_{\alpha, \beta \in \NN_{2r}^{mn}}\widetilde{\iy}_{\alpha+\beta}p_{\alpha}p_{\beta} \\
= \, & \textstyle \sum_{I \subset [m \times n]} \big( \textstyle \sum_{\alpha, \beta \in \NN^{mn}_{r- d_I}}\widetilde{\iy}_{2\gamma_I+2\alpha+2\beta}p_{2\alpha+\gamma_I}p_{2\beta+\gamma_I} \big)\\
= \, & \textstyle \sum_{I \subset [m \times n]} \big( \textstyle \sum_{\alpha, \beta \in \NN^{mn}_{r- d_I}}\iy_{\gamma_I+\alpha+\beta}p_{2\alpha+\gamma_I}p_{2\beta+\gamma_I} \big) = \textstyle \sum_{I \subset [m \times n]}p_I^{\top}\M_{r-d_I}(e_I\iy)p_I \geq 0.
\end{aligned}
\end{equation*}
Here, $p_I$ represents the coefficient vector of the polynomial $p_I(\pi) = \sum_{\alpha \in \NN^{mn}_{r-d_I}}p_{2\alpha+\gamma_I}\pi^{\alpha}$.  The last inequality comes from fact that $\M_{r-d_I}(e_I\iy) \succeq 0$.  Finally, we note that $\iy$ and $\widetilde{\iy}$ evaluate to equal objective values.  By minimizing over $\iy$, we conclude that $\eqref{hierarchy Putinar discrete GW of double degree} \leq \eqref{hierarchy Schmudgen discrete GW}$.

This completes the proof.
\end{proof}

It is worth to observe that from this proof, we can combine both the {\em Schmüdgen-type} hierarchy \eqref{hierarchy Schmudgen discrete GW} and the {\em Putinar-type} hierarchy \eqref{hierarchy Putinar discrete GW} to obtain a new hierarchy, whose number of variables is the same as that of \eqref{hierarchy Schmudgen discrete GW}, but the number of semidefinite constraints is the same as that of \eqref{hierarchy Putinar discrete GW} as following:
\begin{equation} \label{hierarchy Schmudgen-Putinar discrete GW} 
\begin{aligned}
\min \qquad & \sum_{i,k \in [m],\; j,l \in [n]}L_{ij,kl}\cdot \iy_{ij,kl} \\
\mbox{subject to} \qquad & \iy \in \RR^{s(mn,2r)},\;\iy_0 =1, \;\\ 
& \widetilde{\M}_r(\iy) \succeq  0,\\
&\M_{r-1}(r_i(\pi) \iy) = 0, \;\M_{r-1}(r_{i_1}e_{i_2j}(\pi) \iy) =0\; \forall i_1, i_2 \in [m],\ j \in [n]\\
&\M_{r-1}(c_j(\pi) \iy) = 0, \;\M_{r-1}(c_{j_1}e_{ij_2}(\pi) \iy) =0\; \forall j_1,j_2 \in [n],\ i \in [m].
\end{aligned} \tag{SP-DGW-r}
\end{equation}
Here, the matrix $\widetilde{\M}_r(\iy) \in \RR^{s(mn,2r) \times s(mn,2r)}$, where for any $\alpha, \beta \in \NN^{mn}_{2r}$, the $(\alpha,\beta)-$entry is defined by
\begin{equation*}
\widetilde{\M}_r(\iy)(\alpha,\beta) = \iy_{\gamma} \; ~ \text{if}~ \; \exists \gamma\;:\; \alpha+\beta = 2\gamma, \quad \widetilde{\M}_r(\iy)(\alpha,\beta) = 0 ~ \text{otherwise}.
\end{equation*}

\subsection{Dual formulations} 

In this section we describe the dual formulation of \eqref{hierarchy Schmudgen discrete GW} and \eqref{hierarchy Putinar discrete GW}.  

First, we specialize the definition of the preordering to the set $\Pi(\mu_{\CX},\mu_{\CY})$ 
\begin{equation*}
\begin{aligned}
&\CT(\Pi(\mu_{\CX},\mu_{\CY})) = \Big \{\sigma(\pi)+\sum_{I \subset [ m\times n]}e_I\sigma_{I}(\pi) + \sum_{i \in [m]}\lambda_i(\pi) \big (\muxi-\sum_{j \in [n]}\pi_{ij}\big)\\
&+\sum_{j \in [n]}\theta_j(\pi) \big(\muyj-\sum_{i \in [m]}\pi_{ij}\big): \sigma,\sigma_I \in \Sigma[\pi],\; \lambda_i, \theta_j \in \RR[\pi]\; \forall i \in [m],\ j \in [n],\ I \subset [m \times n] \Big \}.
\end{aligned}
\end{equation*}
In a similar fashion, the $r$-th degree truncation is
\begin{equation*}
\begin{aligned}
&\CT(\Pi(\mu_{\CX},\mu_{\CY}))_r =\\
&\Big \{ \sigma(\pi)+\sum_{I \subset [ m\times n]}e_I\sigma_{I}(\pi) + \sum_{i \in [m]}\lambda_i(\pi) \big(\muxi-\sum_{j \in [n]}\pi_{ij}\big)+\sum_{j \in [n]}\theta_j(\pi) \big(\muyj-\sum_{i \in [m]}\pi_{ij} \big):\\
&\sigma\in \Sigma[\pi]_r,\; \sigma_I \in \Sigma[\pi]_{r-\deg e_I},\; \lambda_i, \theta_j \in \RR[\pi]_{r-1}\; \forall i \in [m],\ j \in [n],\ I \subset [m \times n] \Big \}.
\end{aligned}
\end{equation*}

The Schmudgen-type SOS hierarchy is 
\begin{equation*}
\lb(f,\CT(\Pi(\mu_{\CX},\mu_{\CY})))_r = \sup\{c \in \RR:\ f(x) -c \in \CT(\Pi(\mu_{\CX},\mu_{\CY}))_{2r}\}.
\end{equation*}
Note that the SDP instances $\lb(L,\CT(\Pi(\mu_{\CX},\mu_{\CY})))_r$ and \eqref{hierarchy Schmudgen discrete GW} are Lagrangian duals of each other, for every positive integer $r$. 

The Putinar-type moment-SOS hierarchy also lends itself to a dual interpretation.  We specialize the definition of the quadratic module to the set $\widetilde{\Pi}(\mu_{\CX},\mu_{\CY})$
\begin{equation*}
\begin{aligned}
\CQ(\widetilde{\Pi}(\mu_{\CX},\mu_{\CY})) = \Bigl\{\sigma(\widetilde{\pi})+\textstyle\sum_{i \in [m]}\lambda_i(\widetilde{\pi})\big(\muxi-\textstyle\sum_{j \in [n]}\widetilde{\pi}_{ij}^2\big)+\textstyle\sum_{j \in [n]}\theta_j(\widetilde{\pi})\big(\muyj-\textstyle\sum_{i \in [m]}\widetilde{\pi}_{ij}^2\big):\\ \sigma \in \Sigma[\pi],\; \lambda_i,\; \theta_j \in \RR[\widetilde{\pi}]\; \forall i \in [m],\; j \in [n] \Bigr\}.
\end{aligned}
\end{equation*}
In a similar fashion, the $r$-th degree truncation is 
\begin{equation*}
\begin{aligned}
\CQ(\widetilde{\Pi}(\mu_{\CX},\mu_{\CY}))_r = \Bigl\{\sigma(\widetilde{\pi})+\textstyle\sum_{i \in [m]}\lambda_i(\widetilde{\pi})\big(\muxi-\textstyle\sum_{j \in [n]}\widetilde{\pi}_{ij}^2\big)+\textstyle\sum_{j \in [n]}\theta_j(\widetilde{\pi})\big(\muyj-\textstyle\sum_{i \in [m]}\widetilde{\pi}_{ij}^2\big):\\
\sigma \in \Sigma[\pi]_{2r},\; \lambda_i,\; \theta_j \in \RR[\widetilde{\pi}]_{r-2}\; \forall i \in [m],\; j \in [n] \Bigr\}.
\end{aligned}
\end{equation*}
Given any positive integer $r$, the Lagrangian dual of the \eqref{hierarchy Putinar discrete GW} is
\begin{equation*}
    \lb(\widetilde{L},\CQ(\widetilde{\Pi}(\mu_{\CX},\mu_{\CY})))_r = \sup\{c \in \RR: L(x) -c \in \CT(\widetilde{\Pi}(\mu_{\CX},\mu_{\CY})_{2r}\}.
\end{equation*}
A consequence of $\widetilde{\Pi}(\mu_{\CX},\mu_{\CY})$ being Archimedean is that we have strong duality; that is
\begin{equation*}
\lb(\widetilde{L},\CQ(\widetilde{\Pi}(\mu_{\CX},\mu_{\CY})))_r =\eqref{hierarchy Putinar discrete GW} \leq \eqref{eq:gw_discrete}.
\end{equation*}
A consequence of this result is that we also have strong duality for the Schmüdgen-type hierarchy.  

\begin{proposition} \label{thm:strongduality}
Consider the primal-dual instances \eqref{hierarchy Schmudgen discrete GW} and $\lb(L,\CT(\Pi(\mu_{\CX},\mu_{\CY})))_r$.  One has that strong duality holds; that is
\begin{equation*}
\lb(L,\CT(\Pi(\mu_{\CX},\mu_{\CY})))_r = \eqref{hierarchy Schmudgen discrete GW}.
\end{equation*}
\end{proposition}

\noindent {\em Remark: The consequences in Proposition \ref{thm:strongduality} are stated as a standalone result.  We are not aware of more general results that assert strong duality for Schmüdgen-type hierarchies.}

\begin{proof}[Proof of Proposition \ref{thm:strongduality}]
By combining Theorem \ref{connection between 2 hierarchies} and strong duality of the moment-SOS hierarchy under the Archimedean condition (see e.g., \cite{josz2016strong}), we have
\begin{equation*}
\lb(\widetilde{L},\CQ(\widetilde{\Pi}(\mu_{\CX},\mu_{\CY}))_{2r} = \eqref{hierarchy Putinar discrete GW of double degree} = \eqref{hierarchy Schmudgen discrete GW}.
\end{equation*}
As a consequence of weak duality, we also have $\lb(L,\CT(\Pi(\mu_{\CX},\mu_{\CY}))_r \leq \eqref{hierarchy Schmudgen discrete GW}$.  As such, it remains to show that
\begin{equation*}
\lb(L,\CT(\Pi(\mu_{\CX},\mu_{\CY}))_r \geq \lb(\widetilde{L},\CQ(\widetilde{\Pi}(\mu_{\CX},\mu_{\CY}))_{2r}.
\end{equation*}
Let $(c, \sigma, \lambda_i, \theta_j)$ be a feasible solution of $\lb(\widetilde{L},\CQ(\widetilde{\Pi}(\mu_{\CX},\mu_{\CY}))_{2r}$; i.e., 
\begin{equation*}
\widetilde{L}(\widetilde{\pi}) -c = \sigma(\widetilde{\pi})+ \textstyle\sum_{i \in [m]}\lambda_i(\widetilde{\pi})\big(\muxi-\textstyle\sum_{j \in [n]}\widetilde{\pi}_{ij}^2\big)+\textstyle\sum_{j \in [n]}\theta_j(\widetilde{\pi})\big(\muyj-\textstyle\sum_{i \in [m]}\widetilde{\pi}_{ij}^2\big),
\end{equation*}
for some $\sigma \in \Sigma[\pi]_{4r},\; \lambda_i, \theta_j \in \RR[\widetilde{\pi}]_{4r-2}\; \forall i \in [m], j \in [n]$.  In the following, we let $\widetilde{\pi}^2$ denote the vector $(\widetilde{\pi}^2)_{i \in [m], j \in [n]}$.  Also, for any $\varepsilon \in \{1,-1\}^{m \times n}$, we denote $\varepsilon \cdot \widetilde{\pi}= (\varepsilon_{ij}\cdot \widetilde{\pi}_{ij})_{i \in [m], j \in [n]}$.  Using the fact that $\widetilde{L}(\widetilde{\pi}) - c$ is an {\em even} polynomial -- that is, all monomials are of the form $x^{2\alpha}$ -- we have
\begin{equation} \label{eq:squareroot_of_polynomial}
\begin{aligned}
2^{mn}(\widetilde{L}(\widetilde{\pi}) -c) = \textstyle\sum_{\varepsilon \in \{-1,1\}^{m \times n}}\sigma(\varepsilon\cdot \widetilde{\pi})+\Bigl( \textstyle\sum_{\varepsilon \in \{-1,1\}^{m \times n}}\sum_{i \in [m]}\lambda_i(\varepsilon\cdot \widetilde{\pi})\Bigr)\Bigl(\muxi-\textstyle\sum_{j \in [n]}\widetilde{\pi}_{ij}^2\Bigr)\\+\Bigl(\textstyle\sum_{\varepsilon \in \{-1,1\}^{m \times n}}\textstyle\sum_{j \in [n]}\theta_j(\varepsilon\cdot \widetilde{\pi})\Bigr)\Bigl(\muyj-\textstyle\sum_{i \in [m]}\widetilde{\pi}_{ij}^2\Bigr) 
\end{aligned}.
\end{equation}
Since $\lambda_i \in \RR[\widetilde{\pi}]_{4r-2}$ and $\sum_{\varepsilon \in \{-1,1\}^{m \times n}}\sum_{i \in [m]}\lambda_i(\varepsilon\cdot \widetilde{\pi})$ is an even polynomial, there exists a polynomial $\overline{\lambda}_i$ of degree at most $2r-1$ such that 
\begin{equation*}
\textstyle\sum_{\varepsilon \in \{-1,1\}^{m \times n}} \textstyle\sum_{i \in [m]}\lambda_i(\varepsilon\cdot \widetilde{\pi})= \overline{\lambda}_i(\widetilde{\pi}^2) \; \forall i \in [m].
\end{equation*}
Similarly, there exist polynomials $\overline{\theta}_j$ such that 
\begin{equation*}
\textstyle \sum_{\varepsilon \in \{-1,1\}^{m \times n}} \textstyle\sum_{j \in [n]}\theta_j(\varepsilon\cdot \widetilde{\pi}) = \overline{\theta}_j(\widetilde{\pi}^2).
\end{equation*}

We claim that the polynomials $\overline{\lambda}_i$ and $\overline{\theta}_j$ are SOS for all $i \in [m],\ j\in [n]$. The claim can be generalized as the following statement.

{\em Claim:} Let $p \in \Sigma[x]$ with $x \in \RR^n$. Then there exists a polynomial $\overline{p}$ such that
\begin{equation*}
    \overline{p}(x_1^2,\dots,x_n^2) = \textstyle\sum_{\varepsilon \in \{-1,1\}^n} p(\varepsilon \cdot x),\quad \exists p_i \in \Sigma[x]\; \forall i \in [n]\;:\;p(x) = \textstyle \sum_{I \subset [n]}e_I(x)p_I(x).
\end{equation*}
Without loss of generality, it suffices to prove the claim for square polynomial $p^2$ since an SOS polynomial is a sum of squares of polynomials.  We prove this claim by induction on $n$.  We start with the base case $n=1$.  Given a polynomial $p(x)$, separate the monomials to odd and even degrees $p(x) = p_1(x^2) + x p_2(x^2)$.  One then has the identity 
\begin{equation*}
p(x)^2+p(-x)^2 \equiv 2p_1(x^2)^2+ 2x^2p_2(x^2)^2.
\end{equation*}
Therefore, $\overline{p}= 2p_1(x)^2 + 2xp_2(x)^2$ satisfies the claim.  Now assume that the claim holds for $n$.  We need to prove that the claim also holds $n+1$.  We then have 
\begin{equation*}
\begin{aligned}
& p(x_1,\dots,x_{n+1})= p_1(x_1,\dots,x_n,x_{n+1}^2)+x_{n+1}p_2(x_1,\dots,x_n,x_{n+1}^2) \\
\Rightarrow ~& \textstyle\sum_{\varepsilon\in \{-1,1\}}p^2(x_1,\dots,x_n,\varepsilon x_{n+1})= 2p_1^2(x_1,\dots,x_n,x_{n+1}^2)+2x_{n+1}^2p_2^2(x_1,\dots,x_n,x^2_{n+1}) \\
\Rightarrow ~& \textstyle\sum_{\varepsilon\in \{-1,1\}^{n+1}}p^2(\varepsilon \cdot x)= 2\textstyle\sum_{\varepsilon \in \{-1,1\}^n}p_1^2(\varepsilon_1x_1,\dots,\varepsilon_nx_n,x_{n+1}^2) \\
& \qquad \qquad \qquad \qquad \qquad \qquad
 +2x_{n+1}^2 \textstyle\sum_{\varepsilon \in \{-1,1\}^n}p_2^2(\varepsilon_1x_1,\dots,\varepsilon_nx_n,x_{n+1}^2).
\end{aligned}
\end{equation*}
By applying the inductive hypothesis to the last equality we prove the claim.

Subsequently, we conclude that there exists SOS polynomials $\sigma_I$'s with degree at most $2r - d_I$ such that 
\begin{equation*}
\textstyle \sum_{\varepsilon \in \{-1,1\}^{m \times n}}\sigma(\varepsilon\cdot \widetilde{\pi}) = \textstyle \sum_{I \subset [m \times n]}e_I^2(\widetilde{\pi})\sigma_I(\widetilde{\pi}^2).
\end{equation*}
Hence, we have
\begin{equation*}
\widetilde{L}(\widetilde{\pi}) - c = \frac{1}{2^{mn}}\left[\textstyle \sum_{I \subset [m \times n]}e_I^2(\widetilde{\pi})\sigma_I(\widetilde{\pi}^2)+ \textstyle \sum_{i \in [m]}\overline{\lambda}_i(\widetilde{\pi}^2) + \sum_{j \in [n]}\overline{\theta}_j(\widetilde{\pi}^2)\right].
\end{equation*}
We substitute $\widetilde{\pi}^2$ by $\pi$ to obtain 
\begin{equation*}
L(\pi) -c = \frac{1}{2^{mn}}\left[\textstyle \sum_{I \subset [m \times n]}e_I(\pi)\sigma_I(\pi)+\textstyle \sum_{i \in [m]}\overline{\lambda}_i(\pi) + \textstyle \sum_{j \in [n]}\overline{\theta}_j(\pi)\right] \in \CT(\Pi(\mu_{\CX},\mu_{\CY})),
\end{equation*}
which implies that $c \leq \lb(L,\CT(\Pi(\mu_{\CX},\mu_{\CY})))_r$.  Last, we take the supremum over $c$ to conclude
\begin{equation*}
\lb(L,\CT(\Pi(\mu_{\CX},\mu_{\CY}))_r \geq \lb(\widetilde{L},\CQ(\widetilde{\Pi}(\mu_{\CX},\mu_{\CY}))_{2r}. 
\end{equation*}
This completes the proof.
\end{proof}

\section{Convergence Analysis}

The goal of this section is to establish convergence of the proposed SOS hierarchies \eqref{hierarchy Schmudgen discrete GW} and \eqref{hierarchy Putinar discrete GW}.  
The main result of this section is as follows:
\begin{theorem}\label{thm: convergence rate}
Consider the discrete GW problem \eqref{eq:gw_discrete}.  Let $K := \max_{i,k \in [m],\; j,l \in [n]}|L_{ij,kl}|$. Then for any integer $r \geq 4mn$, the lower bound $\lb(L,\CT(\Pi(\mu_{\CX},\mu_{\CY}))_r$ for $\eqref{eq:gw_discrete}$ satisfies:
\begin{equation*}
0 \leq \eqref{eq:gw_discrete} - \lb(L,\CT(\Pi(\mu_{\CX},\mu_{\CY}))_r \leq \frac{C}{r^2} 
(K-\eqref{eq:gw_discrete}+3\sqrt{mn}K\growbeta) + \frac{C^{1/2}}{r} 3\sqrt{mn}K\growbeta .
\end{equation*}
In particular, the hierarchy $\{\lb(L,\CT(\Pi(\mu_{\CX},\mu_{\CY}))_r\}_{r=1}^{\infty}$ converges at the rate of $\mathcal{O}(1/r)$.  
\end{theorem}

Here, $B$ and $C$ are constants which we elaborate on in the following.  The proof of Theorem \ref{thm: convergence rate} builds on ideas in \cite{fang2021sum,slot2111sum,S.o.S-on-simplex, tran2024convergenceratessoshierarchies, tran2025convergence}.  The first step is to apply the following result, which establishes a $\mathcal{O}(1/r^2)$ convergence rate for minimizing functions over the standard simplex via the Schmüdgen-type moment-SOS hierarchy:

\begin{theorem}\label{thm: simplex}
Let $\Delta_n=\{\mathbf{x} \in \mathbb{R}^n:  \sum_{i=1}^nx_i \leq 1, x \geq 0\}$ be the $n$-dimensional standard simplex and let $f \in \mathbb{R}[\mathbf{x}]$ be a polynomial of degree $d$. Then for any $r \geq 2nd$, the lower bound $\lb(f,\mathcal{T}(\mathbf{X}))_r$ for the minimization of $f$ over $\Delta_n$ satisfies:
\begin{equation*}
f_{\min}-\lb(f,\mathcal{T}(\Delta_n))_r \leq \dfrac{C_{\Delta}(n,d)}{r^2}\cdot (f_{\max}-f_{\min}).
\end{equation*}
Here, $f_{\max}$ and $f_{\min}$ denote the maximum and the minimum value of $f$ attained over $\Delta_n$ respectively, while $C_{\Delta}(n, d)$ is a constant depending only on $n, d$.  The constant $C_{\Delta}(n, d)$ can be chosen with polynomial dependence on $n$ for fixed $d$, and with polynomial dependence on $d$ for fixed $n$.
\end{theorem}

The constant $C$ appearing in Theorem \ref{thm: convergence rate} is precisely the constant $C_{\Delta}$ that appears in Theorem \ref{thm: simplex}, and as such we may take $C = C_{\Delta}(mn, 2)$.  

Note that the constraint set in the GW problem \eqref{eq:gw_discrete} is $\Pi(\mu_{\CX},\mu_{\CY})$, which is a strict subset of the probability simplex $\Delta_{mn}$.  As such, we cannot simply apply Theorem \ref{thm: simplex} as it is.  Instead, it is necessary to consider a penalized version of the objective function, defined for values $\lambda > 0$:
\begin{equation}\label{thm: convergence rate inequality 0}
L(\pi,\lambda) := L(\pi) + \lambda \big[ g(\pi)^2 + \textstyle\sum_{i \in [m]}r_i(\pi)^2 + \textstyle\sum_{j \in [n]}c_j(\pi)^2 \big].
\end{equation}

In what follows, we apply Theorem \ref{thm: simplex} with the function $L(\pi,\lambda)$ as our choice of $f$.  Here, $L(\pi,\lambda)$ serves as our proxy for $L(\pi)$ -- note in particular that for any $\pi \in \Pi(\mu_{\CX},\mu_{\CY})$, one has $L(\pi, \lambda) = L(\pi)$.

The introduction of the penalized objective $L(\pi,\lambda)$ incurs an error, which we control using the Łojasiewicz inequality.  In our set-up, the domain of interest $\Pi(\mu_{\CX},\mu_{\CY})$ is {\em polyhedral}, in which case the Łojasiewicz inequality admits a simple form (see e.g., \cite{bergthaller1992distance}).

\begin{lemma}[\cite{bergthaller1992distance},Theorem 0.1]\label{lemma: bound on distance to a polyhedron}
Let $A$ be an $m \times n$ matrix, and let $\beta$ be the smallest number such that for each non-singular sub-matrix $\tilde{A}$ of $A$, one has $|(\tilde{A}^{-1})_{ij}| \leq \beta$; i.e., all entries of $\tilde{A}^{-1}$ are at most $\beta$ in magnitude.  Let $b' \in \RR^m$.  Then for every $b_0 \in \RR^m$ and $x_0 \in \RR^n$ such that 
\begin{equation*}
Ax_0 \leq b_0,
\end{equation*}
there exists $x' \in \RR^n$ satisfying 
\begin{equation*}
Ax' \leq b' \quad \text{and} \quad \|x_0-x'\|_{\infty} \leq n\beta \|b_0 -b'\|_{\infty}.
\end{equation*}
\end{lemma}

Lemma \ref{lemma: bound on distance to a polyhedron} is in effect a special instance of the Łojasiewicz inequality, which is a more general result that provides bounds between a point from the zero locus of a real analytic function.  We apply the result as follows:  First, note that the set $\Pi(\mu_{\CX},\mu_{\CY})$ can be defined as the intersection of these inequalities:
\begin{enumerate}
\item $-\pi_{ij} \leq 0 ~ \forall i \in [m], j \in [n]$,
\item $\sum_{i \in [m]} \pi_{ij} \leq \muyj ~  \forall j \in [n]$,
\item $\sum_{j \in [n]} \pi_{ij} \leq \muxi ~  \forall i \in [m]$,
\item $\sum_{i \in [m], j \in[n]} \pi_{ij} \leq 1$,
\item $- \sum_{i \in [m],j \in[n]} \pi_{ij} \leq -1$.
\end{enumerate}
We express the set $\Pi(\mu_{\CX},\mu_{\CY})$ as a polyhedron in the form $\{ \pi : A (\pi) \leq b' \}$, where $A$ is a linear map and $b'$ is a vector specified via the above inequalities.  Second, let $\pi_{0} \in \Delta_{mn}$ be arbitrary.  Now, $\pi_{0}$ does not necessarily belong to $\Pi(\mu_{\CX},\mu_{\CY})$, but it belongs to the following modified set whereby we replace inequalities (2), (3) and (5) in the above description by these:
\begin{enumerate}
\item[2'.] $\sum_{i \in [m]} \pi_{ij} \leq \sum_{i \in [m]} (\pi_0)_{ij}$,
\item[3'.] $\sum_{j \in [n]} \pi_{ij} \leq \sum_{j \in [n]} (\pi_0)_{ij}$,
\item[5'.] $- \sum_{i \in [m],j \in[n]} \pi_{ij} \leq- \sum_{i \in [m],j \in[n]} (\pi_0)_{ij}$
\end{enumerate}
The resulting set too can be expressed as a polyhedron of the form $\{ \pi : A (\pi) \leq b_0 \}$ with the {\em same} linear map $A$, but with a different vector $b_0$.  Given a closed convex set $S$, let $\mathrm{Proj}(x,S)$ denote the Euclidean projection of $x$ onto $S$.  Then, by using the fact that all norms are equivalent, and by applying Lemma \ref{lemma: bound on distance to a polyhedron}, we have that
\begin{equation} \label{eq:beta_bound}
\begin{aligned}
d(\pi_0, \Pi(\mu_{\CX},\mu_{\CY})) := \| \pi_0 - \mathrm{Proj} (\pi_0, \Pi(\mu_{\CX},\mu_{\CY})) \|  \leq \, & \| \pi_0 - \pi^{\prime} \| \leq \growbeta^{\prime} \| \pi_0 - \pi^{\prime} \|_{\infty} \\
\leq \, & \growbeta \| b^{\prime} - b_0 \|_{\infty} \leq \growbeta h(\pi_0),
\end{aligned}
\end{equation}
for some constants $\growbeta,\growbeta^{\prime}$, and where $\pi^{\prime}$ is the point in $\Pi(\mu_{\CX},\mu_{\CY})$ satisfying the last inequality of Lemma~\ref{lemma: bound on distance to a polyhedron}.  Here, we define
\begin{equation} \label{eq:h_defn}
h(\pi) := \max_{i \in [m],\; j \in [n]} \left\{\;|g(\pi)|,\; |r_i(\pi)|, \; |c_j(\pi)| \right\}.
\end{equation}
(Note in particular that, from the conclusions of Lemma \ref{lemma: bound on distance to a polyhedron}, the constant $\growbeta$ only depends on the linear map $A$ and not the input $\pi_{0}$.)




\begin{proof}[Proof of Theorem \ref{thm: convergence rate}]
We denote $L(\pi,\lambda)_{\min} := \min_{\pi \in \Delta_{mn}} L(\pi,\lambda)$, and $L(\pi,\lambda)_{\max} := \max_{\pi \in \Delta_{mn}} L(\pi,\lambda)$.  Our first steps are to provide bounds on these two quantities.

[Lower bounding $L(\pi,\lambda)_{\min}$]:  The set $\Pi(\mu_{\CX},\mu_{\CY})$ is compact convex.  Denote the projection of $\pi$ to $\Pi(\mu_{\CX},\mu_{\CY})$ by $P(\pi, \Pi(\mu_{\CX},\mu_{\CY}))$.  From \eqref{eq:beta_bound} we have
\begin{equation*}
\|\pi-P(\pi, \Pi(\mu_{\CX},\mu_{\CY}))\| = d(\pi,\Pi(\mu_{\CX},\mu_{\CY})) \leq \growbeta h(\pi).
\end{equation*}

  Based on the definition of $h$ in \eqref{eq:h_defn}, one has $h(\pi)^2 \leq g(\pi)^2 + \sum_{i \in [m]}r_i(\pi)^2 +  \sum_{j \in [n]}c_j(\pi)^2$, and hence
\begin{equation} \label{thm: convergence rate inequality 1}
\begin{aligned}
L(\pi,\lambda) ~\geq~ & \lambda h(\pi)^2 +L(P(\pi,\Pi(\mu_{\CX},\mu_{\CY}))+ \left(L(\pi) - L(P(\pi,\Pi(\mu_{\CX},\mu_{\CY}))\right) \\
~\geq~ & \lambda h(\pi)^2 + \eqref{eq:gw_discrete} - \left|L(\pi) - L(P(\pi,\Pi(\mu_{\CX},\mu_{\CY}))\right|. 
\end{aligned}
\end{equation}
We bound the latter term $|L(\pi) - L(P(\pi,\Pi(\mu_{\CX},\mu_{\CY}))|$.  We do so by providing an upper bound on the Lipschitz constant of $L$ over $\Delta_{mn}$.  Given $\pi' \in \Delta_{mn}$, the gradient of $L$ at $\pi'$ can be bounded as 
\begin{equation*}
\left | \frac{\partial L}{\partial \pi_{ij}}(\pi') \right |= \big|\textstyle\sum_{k \in [m],\; l \in [n]}L_{ij,kl}\pi'_{kl} \big| \leq K \cdot \textstyle\sum_{k \in [m],\; l \in [n]}\pi'_{kl} = K.
\end{equation*}
This implies $\|\frac{\partial L}{\partial \pi}(\pi') \| \leq \sqrt{mn}K$.  Using the inequality $|f(x)-f(y)| \leq \| \nabla f \| \|x-y\|$ on \eqref{thm: convergence rate inequality 1}, we obtain
\begin{equation} \label{thm: convergence rate inequality 2}
\begin{aligned}
L(\pi,\lambda)_{\min} - \eqref{eq:gw_discrete} ~\geq~ & L(\pi,\lambda) - \eqref{eq:gw_discrete} \\
\geq~ & \lambda h(\pi)^2 - \sqrt{mn}Kd(\pi,\Pi(\mu_{\CX},\mu_{\CY})) \\
\geq~ & \lambda h(\pi)^2 - \sqrt{mn}K \growbeta h(\pi) \geq -mnK^2\growbeta^2/(4 \lambda).
\end{aligned}
\end{equation}
The last inequality follows by completing the square.

[Upper bounding $L(\pi,\lambda)_{\max}$]:  We make the following estimates
\begin{enumerate}
\item $L(\pi) = \sum_{i,k \in [m],\; j,l \in [n]}L_{ij,kl}\pi_{ij}\pi_{kl} \leq K ( \sum_{i \in [m],\; j \in [n]}\pi_{ij} )^2 \leq K$,
\item $g(\pi)^2 = (1- \sum_{i \in [m],\; j \in [n]}\pi_{ij})^2 \leq 1$,
\item $\sum_{i\in [m]}r_i(\pi)^2 = \sum_{i \in [m]} (\muxi - \sum_{j \in [n]}\pi_{ij} )^2 \leq \sum_{i \in [m]} \muxi^2 + (\sum_{j \in [n]}\pi_{ij} )^2 \leq 2$, 
\item $\sum_{j\in [n]}c_j(\pi)^2 = \sum_{j \in [n]} (\muyj - \sum_{i \in [m]}\pi_{ij} )^2 \leq \sum_{j \in [n]} \muyj^2 + (\sum_{i \in [m]}\pi_{ij} )^2 \leq 2$.
\end{enumerate}
This gives the $L(\pi,\lambda) \leq K + 5 \lambda$ for all $\pi \in \Delta_{mn}$ and hence
\begin{equation} \label{thm: convergence rate inequality 4}
L(\pi,\lambda)_{\max} \leq K + 5 \lambda.
\end{equation}

[Apply Theorem \ref{thm: simplex}]:  Using the bounds obtained in \eqref{thm: convergence rate inequality 2} and \eqref{thm: convergence rate inequality 4}, we apply Theorem \ref{thm: simplex} with the choice of the quadratic function $L(\pi,\lambda)$ as $f$ to obtain
\begin{equation*}
\begin{aligned}
L(\pi,\lambda)_{\min} - \lb(L(\pi,\lambda), \CT(\Delta_{mn}))_r ~\leq~ & (C_{\Delta}(mn,2)/r^2) (L(\pi,\lambda)_{\max}-L(\pi,\lambda)_{\min}) \\
\leq ~ & (C_{\Delta}(mn,2)/r^2) \big (K + 9\lambda-\eqref{eq:gw_discrete}+ mnK^2\growbeta^2/(4 \lambda) \big ).
\end{aligned}
\end{equation*}
In the following, we let $c_0 := C_{\Delta}(mn,2)/r^2$, and $d_0 := mnK^2\growbeta^2/4$.  By applying the bound \eqref{thm: convergence rate inequality 2} on the LHS, we obtain
\begin{equation*}
\eqref{eq:gw_discrete}-\lb(L(\pi,\lambda), \CT(\Delta_{mn}))_r \, \leq \, c_0  (K-\eqref{eq:gw_discrete} )+5\lambda c_0 + d_0 (c_0+1)/\lambda.
\end{equation*}
Consider the last two terms.  By the Cauchy-Schwarz inequality, there exists a positive number $\lambda^*$ such that 
\begin{equation*} 
5 c_0 \lambda^* + d_0 (c_0+1) / \lambda^* = 2\sqrt{5 c_0 d_0 (c_0+1)} \leq 3\sqrt{mn}K\growbeta(c_0+c_0^2)^{1/2} \leq 3\sqrt{mn}K\growbeta (\sqrt{c_0}+c_0).
\end{equation*}
This gives the inequality
\begin{equation}\label{thm: convergence rate inequality 5}
\eqref{eq:gw_discrete}-\lb(L(\pi,\lambda^*), \CT(\Delta_{mn}))_r \leq c_0 (K-\eqref{eq:gw_discrete} )+3\sqrt{mn}K\growbeta (\sqrt{c_0}+c_0 ).
\end{equation}

[Drawing conclusions]:  For each feasible solution $(c,p) \in \RR \times \CT(\Delta_{mn})_{2r}$ of $\lb(L(\pi,\lambda^*), \CT(\Delta_{mn}))_r$, we can write $L(\pi,\lambda^*)-c = p(\pi)$.  From the definition \eqref{thm: convergence rate inequality 0} we have
\begin{equation*}
L(\pi)-c = p(\pi) -\lambda^* \big [g(\pi)^2 + \textstyle\sum_{i \in [m]}r_i(\pi)^2 +  \textstyle\sum_{j \in [n]}c_j(\pi)^2 \big ] \in \CT(\Pi(\mu_{\CX},\mu_{\CY}))_{2r}.
\end{equation*}
To see why the above expression belongs $\CT(\Pi(\mu_{\CX},\mu_{\CY}))_{2r}$, note that (i) $p \in \CT(\Delta_{mn})_{2r}$, and (ii) one can check via definitions that each of $-g(\pi)^2$, $-\sum_{i \in [m]}r_i(\pi)^2$, $-\sum_{j \in [n]}c_j(\pi)^2$ belong to $\CT(\Delta_{mn})_{2r}$.  
Thus $(c,p(\pi) -\lambda^* [g(\pi)^2 + \sum_{i \in [m]}r_i(\pi)^2 +  \sum_{j \in [n]}c_j(\pi)^2 ]  )$ is a feasible solution of the problem $\lb(L,\CT(\Pi(\mu_{\CX},\mu_{\CY}))_r$, which implies $c \leq \lb(L,\CT(\Pi(\mu_{\CX},\mu_{\CY}))_r$.  By taking the supremum of $c$ over the feasible region of $\lb(L(\pi,\lambda^*), \CT(\Delta_{mn}))_r$, we have
\begin{equation} \label{eq: convergence rate inequality 6}
\lb(L(\pi,\lambda^*), \CT(\Delta_{mn}))_r \leq \lb(L,\CT(\Pi(\mu_{\CX},\mu_{\CY}))_r \leq \eqref{eq:gw_discrete}.
\end{equation}

[Conclusion]:  Finally, by combining \eqref{thm: convergence rate inequality 5} and \eqref{eq: convergence rate inequality 6}, we conclude
\begin{equation*}
0 \leq \eqref{eq:gw_discrete} - \lb(L,\CT(\Pi(\mu_{\CX},\mu_{\CY}))_r \leq c_0  (K-\eqref{eq:gw_discrete} ) + 3\sqrt{mn}K\growbeta (\sqrt{c_0}+c_0).
\end{equation*}
This complete the proof.
\end{proof}

\section{Metric Properties}  \label{sec:metric}

In this section, we investigate metric properties concerning distance measures induced by the SOS hierarchies of the GW problem.  
Given metric measure spaces $\XX=(\CX,d_{\CX},\mu_{\CX})$ and $\YY= (\CY,d_{\CY},\mu_{\CY})$, recall the $L_p$ distortion distance $\Delta_p$ given by
\begin{equation*}
\Delta_p(\XX,\YY) = \Big (\min_{\pi \in \Pi(\mu_{\CX},\mu_{\CY})}\int_{(\CX \times \CY)^2}\left|d_{\CX}(x,x')-d_{\CY}(y,y') \right|^pd\pi \otimes \pi \Big )^{1/p}.
\end{equation*}
Our earlier discussion on the SOS hierarchies of the GW problem suggest natural analogs of the $L_p$ distance that are tractable to compute.  More precisely, we define
\begin{equation*}
\Delta_{p,r} := \eqref{hierarchy Schmudgen discrete GW}(\XX,\YY)^{1/p} \qquad \text{where} \qquad L(x,x',y,y') = |d_{\CX}(x,x')-d_{\CY}(y,y')|^p.
\end{equation*}

The main result of this section is to show that $\Delta_{p,r}$ defines a {\em pseudo-metric}; that is, it satisfies
\begin{itemize}
    \item $\Delta_{p,r}(\XX,\XX) =0$,
    \item (Non-negativity) $\Delta_{p,r}(\XX,\YY) \geq 0 $ for all $\XX$, $\YY$,
    \item (Symmetry) $\Delta_{p,r}(\XX,\YY)= \Delta_{p,r}(\YY,\XX)$, and
    \item (Triangle Inequality) $\Delta_{p,r}(\XX,\ZZ) \leq \Delta_{p,r}(\XX,\YY) + \Delta_{p,r}(\YY,\ZZ)$.
\end{itemize}

Among these properties, non-negativity and symmetry are straightforward to prove and we assert these without proof.  The triangle inequality is more troublesome, and is the focus of this section.  The key pre-cursor to proving the triangle inequality is Lemma \ref{thm:gluinglemma}, which best understood as `gluing lemma'-type of result that explains how pairs of moment sequences can be `glued' together to form valid moment sequences.  Lemma \ref{thm:gluinglemma} is the analog of the so-called `gluing lemma' that appears throughout optimal transport; see, for instance \cite[Lemma 7.6]{villani2003topics} in the context of classical optimal transport problems, as well as \cite[Lemma 1.4]{Sturm:23} in the context of the Gromov-Wasserstein problem.  The gluing lemma is a basic result in optimal transport that underlies many useful properties, such as triangle inequalities and geodesics of the metric-measure space.  

We begin by making a series of simplifying observations.  Consider metric measure spaces $\XX$ and $\YY$.  To simplify notation, we denote $\alpha := \mu_{\CX}$ and $\beta := \mu_{\CY}$.  We assume that these distributions have finite support $\{x_1,\dots,x_m\}$ and $\{y_1,\dots,y_n\}$ respectively.  As such, the metric measure spaces in question are $\XX= (\CX,d_{\CX},\alpha)$ and $\YY= (\CY,d_{\CY},\beta)$ respectively. Recall the feasible set of \eqref{hierarchy Schmudgen discrete GW}

\begin{multline*}
    \Pi_r(\alpha,\beta) := \Bigl\{ \iy \in \RR^{s(mn,2r)}:\;\iy_0 =1,\; \M_{r-d_I}(e_I\iy) \succeq 0 \; \forall I \subset [m \times n],\ \M_{r-1}(r_i\iy) =0, \\ \M_{r-1}(c_j\iy) =0,\;
\; \M_{r-1}(r_{i_1}e_{i_2j}\iy) =0 ,\ \M_{r-1}(c_{j_1}e_{ij_2}\iy)=0 \;\forall i_1,i_2,i \in [m],\ j,j_1,j_2\in [n]\Bigr\}.
\end{multline*}

We identify the components of $\iy$ via the Riesz functional $\ell_{\iy}: \; \RR[\pi] \to \RR$ defined by 
\begin{equation*}
    \ell_{\iy}(\pi^{\gamma}) = \iy_{\gamma}, \; \forall \gamma \in \NN^{mn}.
\end{equation*}
Then the constraints defining $\Pi_r(\alpha,\beta)$ can be reduced as follows:   First, recall the linear constraints stemming from the localizing matrices $\M_{r-1}(r_i\iy) =0, \; \M_{r-1}(c_j\iy) =0,\;
\; \M_{r-1}(r_{i_1}e_{i_2j}\iy) =0 ,\ \M_{r-1}(c_{j_1}e_{ij_2}\iy)=0 \;\forall i_1,i_2,i \in [m],\ j,j_1,j_2\in [n]$.  By summing these constraints over $i$, we have $\M_{r-1}((1 - \sum_{i \in [m],j \in [n]}\pi_{ij})e_{i^{\prime} j^{\prime}}\iy) = 0 \; \forall i^{\prime} \in [m],\; j^{\prime} \in [m]$.  This is equivalent to saying that 
\begin{equation*}
\ell_{\iy}(\pi^{\tau}) = \sum_{i \in [m], j\in [n]} \ell_{\iy}(\pi^{\tau +e_{ij}}) \quad \forall \tau \in \NN^{mn}_{2r-1},
\end{equation*}
where $\{e_{ij}\}_{i \in [m], j\in [n]}$ is the standard basis of $\RR^{mn}$.  By iterating to sums of more than one term of the form $e_{ij}$, we obtain the following relationship: for all $\tau \in \NN^{mn}_{2r}$ such that $|\tau| = 2r-l$, one has
\begin{equation*}
\ell_{\iy}(\pi^{\tau}) = \sum_{\substack{i_1,\dots, i_l \in [m]  \\ j_1,\dots,j_l \in [n]}}\ell_{\iy}(\pi^{\tau +\sum_{s=1}^le_{i_sj_s}}).
\end{equation*}
There are repeated terms in the RHS, and so we group the sum by distinct terms.  Concretely, given $\tau = (\tau_1,\dots,\tau_{mn}) \in \NN^{mn}_r$ of length $t$, we define $c(\tau) := \binom{l!}{\tau_1!\cdots \tau_n!}$ as the number of tuples $(e_{i_1j_1},\dots,e_{i_{t}j_{t}})$ such that $\sum_{s=1}^{t}e_{i_sj_s}= \tau$.  Define $\overline{\NN}^n_r= \{\tau \in \NN^{n}_r \;:\; |\tau| =r\}$. Then for any pair of integers $(l_1,l_2)$ such that $l_1+l_2=l$, we can rewrite the component of $\iy$ indexed by $\tau$ of the length $|\tau|= 2r-l$ as 
\begin{equation}\label{decomposition of moment sequence}
    \ell_{\iy}(\pi^{\tau})= \sum_{\sigma \in \overline{\NN}^n_{l_1},\theta \in \overline{\NN}^n_{l_2}}c(\sigma)c(\theta)\ell_{\iy}(\pi^{\tau+\sigma+\theta}).
\end{equation}
The following result describes how we can reduce the size of the moment matrix and the localizing matrices.
\begin{proposition}\label{prop reduciton on semidefinite constrains}
Let $\iy \in \RR^{s(mn,2r)}$ be an element of $\Pi_r(\alpha,\beta)$; i.e., the following conditions are satisfied: (i) $\iy_0 =1$, (ii) $\M_{r-1}(r_{i_1}e_{i_2j}(\pi) \iy) = 0 \; \forall i_1,i_2 \in [m],\; j \in [n]$, and (iii) $\M_{r-1}(c_{j_1}e_{ij_2}(\pi) \iy) = 0 \; \forall i \in [m],\; j_1, j_2 \in [n]$.  Let $\overline{\M}_{r-d_I}(e_I\iy)$ denote the principal sub-matrix of $\M_{r-d_I}(e_I\iy)$ obtained by removing all columns and rows indexed by $\tau \in \NN^{mn}_{r-d_I-1}$.  One then has
$$
\M_{r - d_I}(e_I\iy) \succeq 0 \; \forall I \subset [m \times n] \quad \Leftrightarrow \quad \overline{\M}_{r-d_I}(e_I\iy) \succeq 0 \; \forall I \subset [m \times n].
$$
\end{proposition}

\begin{proof}
First note that $\overline{\M}_{r-d_I}(e_I\iy)$ is a principal sub-matrix of $\M_{r-d_I}(e_I\iy)$.  Hence if $\M_{r-d_I}(e_I\iy) \succeq 0$, we immediately have that $\overline{\M}_{r-d_I}(e_I\iy) \succeq 0$.  As such, it suffices to simply prove the reverse implication $
\overline{\M}_{r-d_I}(e_I\iy) \succeq 0 \; \Rightarrow\; \M_{r-d_I}(e_I\iy) \succeq 0$.

Suppose $\overline{\M}_{r-d_I}(e_I\iy) \succeq 0$.  For any $\vv=(v_{\tau})_{\tau \in \NN^{mn}_{r-d_I}} \in \RR^{s(mn,r-d_I)}$, we have the following
\begin{equation*}
\begin{aligned}
\vv^{\top} \M_{r-d_I}(e_I\iy)\vv \overset{(a)}{=} & \sum_{\tau, \sigma \in \NN^{mn}_{r-d_I}}v_{\tau}v_{\sigma}\ell_{\iy}(e_I\pi^{\tau +\sigma})\\
=\;& \sum_{\tau, \sigma \in \NN^{mn}_{r-d_I}}v_{\tau}v_{\sigma}\sum_{\substack{\theta \in \overline{\NN}^n_{r-d_I-|\tau|}\\ \zeta \in \overline{\NN}^n_{r-d_I-|\sigma|}}}c(\theta)c(\zeta)\ell_{\iy}(e_I\pi^{(\tau +\theta)+(\sigma+\zeta)})\\
=\;& \sum_{\tau, \sigma \in \overline{\NN}^n_{r-d_I}}\sum_{\theta \preceq \tau,\;  \zeta \preceq \sigma}[c(\tau - \theta)v_{\theta}][c(\sigma - \zeta)v_{\zeta}]\ell_{\iy}(e_I\pi^{\tau + \sigma})\\
=\;& \sum_{\tau, \sigma \in \overline{\NN}^n_{r-d_I}}\Big(\sum_{\theta \preceq \tau}c(\tau - \theta)v_{\theta}\Big)\Big(\sum_{\zeta \preceq \sigma}c(\sigma - \zeta)v_{\zeta}\Big)\ell_{\iy}(e_I\pi^{\tau + \sigma}) \overset{(b)}{\geq} 0.
\end{aligned}
\end{equation*}
Here, we use the notation $\theta \preceq \tau$ to imply component-wise inequalities; i.e., $\theta_k \leq \tau_k\; \forall k$, where $\theta_k$ and $\tau_k$ denote the corresponding components of $\theta$ and $\tau$. (a) follows from \eqref{decomposition of moment sequence}, and (b) follows from the assumption that $\overline{\M}_{r-d_I}(e_I\iy) \succeq 0$.  Since $\vv$ is arbitrary, we conclude that $\M_{r-d_I} \succeq 0$.
\end{proof}

It remains to prove the triangular inequality.  In the remainder of this section we adopt the following notation.  We consider metric measure spaces $\XX= (\CX,d_{\CX},\alpha)$ with $\CX = \{ x_1,\dots, x_m\}$, $\YY = (\CY,d_{\CY},\beta)$ with $\CY = \{ y_1,\dots,y_n\}$, and $\ZZ = (\CZ,d_{\CZ},\gamma)$ with $\CZ = \{ z_1,\dots,z_p\}$. We let $\pi^{(1)}:= (\pi^{(1)}_{ij})_{i \in [m], j\in [n]}$, $\pi^{(2)}:= (\pi^{(2)}_{jk})_{i \in [n], j\in [p]}$, $\pi^{(3)}:= (\pi^{(3)}_{ki})_{k\in [p], i \in [m]}$, and $\pi := (\pi_{ijk})_{i \in [m], j\in [n], k \in [p]}$ be variables in $\RR^{m \times n}$, $\RR^{n \times p}$, $\RR^{p \times m}$, and $\RR^{m \times n \times p}$, respectively.  In what follows, we use the Riesz functional to identify components of moment sequences.

\begin{lemma}[Gluing lemma] \label{thm:gluinglemma} Let $\iy_1 \in \Pi_r(\alpha,\beta)$ and $\iy_2 \in \Pi_r(\beta,\gamma)$. Then there exists $\iy \in \RR^{s(mnp,2r)}$ satisfying the following for all $1 \leq l \leq 2r$: 
\begin{equation*}
\ell_{\iy^{(1)}} \Big (\prod_{s=1}^l\pi^{(1)}_{i_sj_s} \Big )= \sum_{k_1,\dots, k_l \in [p]}\ell_{\iy} \Big ( \prod_{s=1}^l\pi_{i_sj_sk_s} \Big ), \quad \ell_{\iy^{(2)}} \Big (\prod_{s=1}^l\pi^{(1)}_{j_sk_s} \Big )= \sum_{i_1,\dots, i_l \in [m]}\ell_{\iy} \Big ( \prod_{s=1}^l\pi_{i_sj_sk_s} \Big ).
\end{equation*}
In addition, the sequence $\iy^{(3)} \in \RR^{s(pm,2r)}$, defined by 
\begin{equation*}
\ell_{\iy^{(3)}} \Big (\prod_{s=1}^l\pi^{(3)}_{k_si_s} \Big )= \sum_{j_1,\dots, j_l \in [n]}\ell_{\iy} \Big ( \prod_{s=1}^l\pi_{i_sj_sk_s} \Big ),
\end{equation*}
belongs to $\Pi_r(\gamma,\alpha)$.
\end{lemma}

\begin{proof}[Proof of Lemma \ref{thm:gluinglemma}]
We begin defining $\iy \in \RR^{s(mnp,2r)}$ as the moment sequence satisfying for any $1 \leq l \leq 2r$, we observe the following for any $1 \leq l \leq 2r$
\begin{equation*}
\ell_{\iy} \Big ( \prod_{s=1}^l\pi_{i_sj_sk_s} \Big )= \begin{cases}
    0 & \text{if }\prod_{s=1}^l\beta_{j_s}=0,\\
    \frac{\ell_{\iy^{(1)}}\Big ( \prod_{s=1}^l\pi^{(1)}_{i_sj_s} \Big )\ell_{\iy^{(2)}} \Big ( \prod_{s=1}^l\pi^{(2)}_{j_sk_s} \Big )}{\prod_{s=1}^l\beta_{j_s}} & \text{otherwise.}
\end{cases} 
\end{equation*}
Without loss of generality, we may omit the atoms of zero mass out of the marginal spaces.  As such, for simplicity, we simply assume that $\beta_j \neq 0$ for all $j \in [n]$.  In what follows, we always use the latter definition for $\ell_{\iy} \Big ( \prod_{s=1}^l\pi_{i_sj_sk_s} \Big )$.

First we note that
\begin{align*}
        & \sum_{k_1,\dots, k_l \in [p]}\ell_{\iy}\left( \prod_{s=1}^l\pi_{i_sj_sk_s}\right)=\sum_{k_1,\dots, k_l \in [p]} \frac{\ell_{\iy^{(1)}}\left( \prod_{s=1}^l\pi^{(1)}_{i_sj_s}\right)\ell_{\iy^{(2)}}\left( \prod_{s=1}^l\pi^{(2)}_{j_sk_s}\right)}{\prod_{s=1}^l\beta_{j_s}}\\
        =\; & \ell_{\iy^{(1)}}\left( \prod_{s=1}^l\pi^{(1)}_{i_sj_s}\right)\left(\sum_{k_1,\dots, k_l \in [p]} \frac{\ell_{\iy^{(2)}}\left( \prod_{s=1}^l\pi^{(2)}_{j_sk_s}\right)}{\prod_{s=1}^l\beta_{j_s}}\right)\\
        =\;& \ell_{\iy^{(1)}}\left( \prod_{s=1}^l\pi^{(1)}_{i_sj_s}\right)\left(\sum_{k_1,\dots, k_{l-1} \in [p]} \frac{\sum_{k_l \in [p]}\ell_{\iy^{(2)}}\left( \prod_{s=1}^l\pi^{(2)}_{j_sk_s}\right)}{\prod_{s=1}^l\beta_{j_s}}\right)\\
         =\;& \ell_{\iy^{(1)}}\left( \prod_{s=1}^l\pi^{(1)}_{i_sj_s}\right)\left(\sum_{k_1,\dots, k_{l-1} \in [p]} \frac{\beta_{l}\ell_{\iy^{(2)}}\left( \prod_{s=1}^{l-1}\pi^{(2)}_{j_sk_s}\right)}{\prod_{s=1}^l\beta_{j_s}}\right)\\
         =\;& \ell_{\iy^{(1)}}\left( \prod_{s=1}^l\pi^{(1)}_{i_sj_s}\right)\left(\sum_{k_1,\dots, k_{l-1} \in [p]} \frac{\ell_{\iy^{(2)}}\left( \prod_{s=1}^{l-1}\pi^{(2)}_{j_sk_s}\right)}{\prod_{s=1}^{l-1}\beta_{j_s}}\right) = \dots =\; \ell_{\iy^{(1)}}\left( \prod_{s=1}^l\pi^{(1)}_{i_sj_s}\right).
\end{align*}
In a similar fashion, we also have
\begin{equation*}
\ell_{\iy^{(2)}}\left(\prod_{s=1}^l\pi^{(1)}_{j_sk_s} \right)= \sum_{i_1,\dots, i_l \in [m]}\ell_{\iy}\left( \prod_{s=1}^l\pi_{i_sj_sk_s}\right).
\end{equation*}

Next we proceed to check that $\iy^{(3)} \in \Pi_r(\gamma,\alpha)$.  First we have
    \begin{align*}
         &\sum_{i_s \in [m]}\ell_{\iy^{(3)}}\left(\prod_{s=1}^l\pi^{(3)}_{k_si_s} \right)=\sum_{i_s \in [m]}\sum_{j_1,\dots, j_l \in [n]}\ell_{\iy}\left( \prod_{s=1}^l\pi_{i_sj_sk_s}\right)\\
         =\;&\sum_{i_s \in [m]}\sum_{j_1,\dots, j_l \in [n]}\frac{\ell_{\iy^{(1)}}\left( \prod_{s=1}^l\pi^{(1)}_{i_sj_s}\right)\ell_{\iy^{(2)}}\left( \prod_{s=1}^l\pi^{(2)}_{j_sk_s}\right)}{\prod_{s=1}^l\beta_{j_s}}\\
         =\;&\sum_{j_1,\dots, j_{l}}\frac{\left(\sum_{i_s \in [m]}\ell_{\iy^{(1)}}\left( \prod_{s=1}^l\pi^{(1)}_{i_sj_s}\right)\right)\ell_{\iy^{(2)}}\left( \prod_{s=1}^l\pi^{(2)}_{j_sk_s}\right)}{\prod_{s=1}^l\beta_{j_s}}\\
         =\;&\sum_{j_1,\dots, j_{l}}\frac{\alpha_l\ell_{\iy^{(1)}}\left( \prod_{s=1}^l\pi^{(1)}_{i_sj_s}\right)\ell_{\iy^{(2)}}\left( \prod_{s=1}^l\pi^{(2)}_{j_sk_s}\right)}{\prod_{s=1}^l\beta_{j_s}}\\
          =\;&\alpha_l\sum_{j_1,\dots, j_{l-1} \in [n]}\frac{\ell_{\iy^{(1)}}\left( \prod_{s=1}^{l-1}\pi^{(1)}_{i_sj_s}\right)\left(\sum_{j_l \in [n]}\ell_{\iy^{(2)}}\left( \prod_{s=1}^l\pi^{(2)}_{j_sk_s}\right)\right)}{\prod_{s=1}^l\beta_{j_s}}\\
          =\;&\alpha_l\sum_{j_1,\dots, j_{l-1} \in [n]}\frac{\ell_{\iy^{(1)}}\left( \prod_{s=1}^l\pi^{(1)}_{i_sj_s}\right)\ell_{\iy^{(2)}}\left( \prod_{s=1}^{l-1}\pi^{(2)}_{j_sk_s}\right)}{\prod_{s=1}^{l-1}\beta_{j_s}}= \alpha_l \ell_{\iy^{(3)}}\left(\prod_{s=1}^{l-1}\pi^{(3)}_{k_si_s} \right).
    \end{align*}
    This implies that $\M_{r-1}\big(\big(\sum_{k \in [p]}\pi^{(3)}_{ki}-\alpha_i\big)\iy^{(3)}\big)=0$ and $\M_{r-1}\big(\big(\sum_{k_1 \in [p]}\pi^{(3)}_{k_1i_1}-\alpha_{i_1}\big)e_{k_2i_2}\iy^{(3)}\big)=0$ for all $i,i_1,i_2 \in [m],\; k_2\in [p]$.  In a similar fashion we also have:
\begin{equation*}
\M_{r-1}\Big(\Big(\sum_{i \in [m]}\pi^{(3)}_{ki}-\gamma_k\Big)\iy^{(3)}\Big)=0, \M_{r-1}\Big(\Big(\sum_{i_1 \in [m]}\pi^{(3)}_{k_1i_1}-\gamma_{k_1}\Big)e_{k_2i_2}\iy^{(3)}\Big)=0 \; \forall i_2 \in [m], k,k_1,k_2 \in [p].
\end{equation*}

It remains to prove the semidefinite constraints on $\iy^{(3)}$.  This step is actually not as difficult as the following sequence of steps may appear to suggest.  The key idea we need is the basic fact that every PSD matrix admits a (unique) PSD matrix square root. 

Based on Proposition \eqref{prop reduciton on semidefinite constrains}, it suffices to check that $\overline{\M}_{r-d_K}(e^{(3)}_K\iy^{(3)}) \succeq 0$ for all $K \subset [p \times m]$.  To simplify notation, we fix $r-d_K=t$ and $\deg e_K^{(3)} = d$.  We then define $K:= \{(k^*_s,i^*_s) \in [p \times m]:\ s \in [d] \}$ so that we can write $e_K^{(3)}= \prod_{s=1}^d \pi_{k_s^*i_s^*}$.  For all vectors $\vv$, we use $v_{\tau}$ or $v(\pi^{\tau})$ to denote the component of $\vv$ that is indexed by $\tau$.  Note that
\begin{equation*}
\begin{aligned}
& \vv^{\top}\overline{\M}_{t}(e^{(3)}_K\iy^{(3)})\vv = \sum_{\tau, \sigma \in \overline{\NN}^{pm}_t}v_{\tau}v_{\sigma}\ell_{\iy^{(3)}}(e^{(3)}_K\pi^{\tau +\sigma})\\
=\;& \sum_{\substack{k_1,\dots,k_t \in[p],\;i_1,\dots,i_t \in [m]\\k'_1,\dots,k'_t \in[p],\;i'_1,\dots,i'_t \in [m]}}\frac{v(\prod_{s=1}^t\pi^{(3)}_{k_si_s})}{c(\prod_{s=1}^t\pi^{(3)}_{k_si_s})}\frac{v(\prod_{s=1}^t\pi^{(3)}_{k'_si'_s})}{c(\prod_{s=1}^t\pi^{(3)}_{k'_si'_s})}\ell_{\iy^{(3)}}\left(\prod_{s=1}^d \pi_{k_s^*i_s^*}\prod_{s=1}^t\pi^{(3)}_{k_si_s}\prod_{s=1}^t\pi^{(3)}_{k'_si'_s}\right).
\end{aligned}
\end{equation*}

Denote $\frac{v(\prod_{s=1}^t\pi^{(3)}_{k_si_s})}{c(\prod_{s=1}^t\pi^{(3)}_{k_si_s})}=:v'(\prod_{s=1}^t\pi^{(3)}_{k_si_s})$ so that the above expression simplifies to
\begin{align*}
& \vv^{\top}\overline{\M}_{t}(e^{(3)}_K\iy^{(3)})\vv\\
=\;& \sum_{\substack{k_1,\dots,k_t \in[p],\;i_1,\dots,i_t \in [m]\\k'_1,\dots,k'_t \in[p],\;i'_1,\dots,i'_t \in [m]}}v'(\prod_{s=1}^t\pi^{(3)}_{k_si_s})v'(\prod_{s=1}^t\pi^{(3)}_{k'_si'_s})\ell_{\iy^{(3)}}\left(\prod_{s=1}^d \pi_{k_s^*i_s^*}\prod_{s=1}^t\pi^{(3)}_{k_si_s}\prod_{s=1}^t\pi^{(3)}_{k'_si'_s}\right)\\
=\;& \sum_{\substack{j_1^*,\dots,j_d^* \in [n]\\ j_1,\dots,j_t \in [n]\\ j'_1,\dots,j'_t \in [n]}}\sum_{\substack{k_1,\dots,k_t \in[p]\\i_1,\dots,i_t \in [m]\\k'_1,\dots,k'_t \in[p]\\i'_1,\dots,i'_t \in [m]}}v'(\prod_{s=1}^t\pi^{(3)}_{k_si_s})v'(\prod_{s=1}^t\pi^{(3)}_{k'_si'_s})\ell_{\iy^{(3)}}\left(\prod_{s=1}^d \pi_{k_s^*j^*_si_s^*}\prod_{s=1}^t\pi_{k_sj_si_s}\prod_{s=1}^t\pi_{k'_sj'_si'_s}\right)\\
=\;& \sum_{\substack{j_1^*,\dots,j_d^* \in [n]\\ j_1,\dots,j_t \in [n]\\ j'_1,\dots,j'_t \in [n]}}\sum_{\substack{k_1,\dots,k_t \in[p]\\i_1,\dots,i_t \in [m]\\k'_1,\dots,k'_t \in[p]\\i'_1,\dots,i'_t \in [m]}}v'(\prod_{s=1}^t\pi^{(3)}_{k_si_s})v'(\prod_{s=1}^t\pi^{(3)}_{k'_si'_s})\\
& \quad \times \frac{\ell_{\iy^{(1)}}\left( \prod_{s=1}^d \pi^{(1)}_{i^*_sj_s^*}\prod_{s=1}^t\pi^{(1)}_{i_sj_s}\prod_{s=1}^t\pi^{(1)}_{i'_sj'_s}\right)\ell_{\iy^{(2)}}\left( \prod_{s=1}^d \pi^{(2)}_{j_s^*k^*_s}\prod_{s=1}^t\pi^{(2)}_{j_sk_s}\prod_{s=1}^t\pi^{(2)}_{j'_sk'_s}\right)}{\prod_{s=1}^d\beta_{j^*_s}\prod_{s=1}^t\beta_{j_s}\prod_{s=1}^t\beta_{j'_s}}\\
=\;& \sum_{\substack{j_1^*,\dots,j_d^* \in [n]\\ j_1,\dots,j_t \in [n],\; j'_1,\dots,j'_t \in [n]\\i_1,\dots,i_t \in [m],\;i'_1,\dots,i'_t \in [m]}}\frac{\ell_{\iy^{(1)}}\left( \prod_{s=1}^d \pi^{(1)}_{i^*_sj_s^*}\prod_{s=1}^t\pi^{(1)}_{i_sj_s}\prod_{s=1}^t\pi^{(1)}_{i'_sj'_s}\right)}{\prod_{s=1}^d\beta_{j^*_s}\prod_{s=1}^t\beta_{j_s}\prod_{s=1}^t\beta_{j'_s}}\\
& \quad \times \sum_{\substack{k_1,\dots,k_t \in[p]\\k'_1,\dots,k'_t \in[p]}}v'(\prod_{s=1}^t\pi^{(3)}_{k_si_s})v'(\prod_{s=1}^t\pi^{(3)}_{k'_si'_s}) \ell_{\iy^{(2)}}\left( \prod_{s=1}^d \pi^{(2)}_{j_s^*k^*_s}\prod_{s=1}^t\pi^{(2)}_{j_sk_s}\prod_{s=1}^t\pi^{(2)}_{j'_sk'_s}\right) .
\end{align*}

Let $J:= \{ (j^*_1,k^*_1),\dotsm,(j^*_d,k^*_d)\}$.  Since $\overline{\M}_t(e^{(2)}_J\iy^{(2)}) \succeq 0$, there exists $Q \succeq 0$ such that $\overline{\M}_t(e^{(2)}_J\iy^{(2)})= Q^2$.  We write the components of $\overline{\M}_t(e^{(2)}_J\iy^{(2)})$ as 
\begin{align*}
    &\ell_{\iy^{(2)}}\left( \prod_{s=1}^d \pi^{(2)}_{j_s^*k^*_s}\prod_{s=1}^t\pi^{(2)}_{j_sk_s}\prod_{s=1}^t\pi^{(2)}_{j'_sk'_s}\right)= \sum_{\tau \in \overline{\NN}^{np}_t}Q\left(\prod_{s=1}^t\pi^{(2)}_{j_sk_s},\pi^{(2)\tau}\right)Q\left(\prod_{s=1}^t\pi^{(2)}_{j'_sk'_s},\pi^{(2)\tau}\right)\\
    =\;& \sum_{\substack{j''_1,\dots,j''_t \in [n]\\ k''_1,\dots,k''_t \in [p]}}\frac{1}{c(\prod_{s=1}^t\pi^{(2)}_{j''_sk''_s})}Q\left(\prod_{s=1}^t\pi^{(2)}_{j_sk_s},\prod_{s=1}^t\pi^{(2)}_{j''_sk''_s}\right)Q\left(\prod_{s=1}^t\pi^{(2)}_{j'_sk'_s},\prod_{s=1}^t\pi^{(2)}_{j''_sk''_s}\right)\\
    =\;& \sum_{\substack{j''_1,\dots,j''_t \in [n]\\ k''_1,\dots,k''_t \in [p]}}Q'\left(\prod_{s=1}^t\pi^{(2)}_{j_sk_s},\prod_{s=1}^t\pi^{(2)}_{j''_sk''_s}\right)Q'\left(\prod_{s=1}^t\pi^{(2)}_{j'_sk'_s},\prod_{s=1}^t\pi^{(2)}_{j''_sk''_s}\right),
\end{align*}
whereby we define $Q' \big (\prod_{s=1}^t\pi^{(2)}_{j_sk_s},\prod_{s=1}^t\pi^{(2)}_{j''_sk''_s} \big)=\frac{Q\left(\prod_{s=1}^t\pi^{(2)}_{j_sk_s},\prod_{s=1}^t\pi^{(2)}_{j''_sk''_s}\right)}{\sqrt{c(\prod_{s=1}^t\pi^{(2)}_{j''_sk''_s})}}$.  We combine with the above equality to obtain 
\begin{align*}
    &\sum_{\substack{k_1,\dots,k_t \in[p]\\k'_1,\dots,k'_t \in[p]}}v'(\prod_{s=1}^t\pi^{(3)}_{k_si_s})v'(\prod_{s=1}^t\pi^{(3)}_{k'_si'_s}) \ell_{\iy^{(2)}}\left( \prod_{s=1}^d \pi^{(2)}_{j_s^*k^*_s}\prod_{s=1}^t\pi^{(2)}_{j_sk_s}\prod_{s=1}^t\pi^{(2)}_{j'_sk'_s}\right)\\
    =\;& \sum_{\substack{k_1,\dots,k_t \in[p]\\k'_1,\dots,k'_t \in[p]\\j''_1,\dots,j''_t \in [n]\\ k''_1,\dots,k''_t \in [p]}}v'(\prod_{s=1}^t\pi^{(3)}_{k_si_s})v'(\prod_{s=1}^t\pi^{(3)}_{k'_si'_s})Q'\left(\prod_{s=1}^t\pi^{(2)}_{j_sk_s},\prod_{s=1}^t\pi^{(2)}_{j''_sk''_s}\right)Q'\left(\prod_{s=1}^t\pi^{(2)}_{j'_sk'_s},\prod_{s=1}^t\pi^{(2)}_{j''_sk''_s}\right)\\
    =\;& \sum_{\substack{j''_1,\dots,j''_t \in [n]\\ k''_1,\dots,k''_t \in [p]}}\left[\sum_{k_1,\dots,k_t \in [m]}v'(\prod_{s=1}^t\pi^{(3)}_{k_si_s})Q'\left(\prod_{s=1}^t\pi^{(2)}_{j_sk_s},\prod_{s=1}^t\pi^{(2)}_{j''_sk''_s}\right)\right]\\
    \cdot& \left[\sum_{k'_1,\dots,k'_t \in [m]}v'(\prod_{s=1}^t\pi^{(3)}_{k'_si'_s})Q'\left(\prod_{s=1}^t\pi^{(2)}_{j'_sk'_s},\prod_{s=1}^t\pi^{(2)}_{j''_sk''_s}\right)\right]\\
    =\;& \sum_{\substack{j''_1,\dots,j''_t \in [n]\\ k''_1,\dots,k''_t \in [p]}}g(\prod_{s=1}^t\pi^{(2)}_{j''_sk''_s},\prod_{s=1}^t\pi^{(1)}_{i_sj_s})g(\prod_{s=1}^t\pi^{(2)}_{j''_sk''_s},\prod_{s=1}^t\pi^{(1)}_{i'_sj'_s}),
\end{align*}
where $g(\prod_{s=1}^t\pi^{(2)}_{j''_sk''_s},\prod_{s=1}^t\pi^{(1)_{i_sj_s}}) := \sum_{k_1,\dots,k_t \in [m]}v'(\prod_{s=1}^t\pi^{(3)}_{k_si_s})Q' \big (\prod_{s=1}^t\pi^{(2)}_{j_sk_s},\prod_{s=1}^t\pi^{(2)}_{j''_sk''_s} \big )$. Hence, we obtain that 
\begin{align*}
    & \vv^{\top}\overline{\M}_{t}(e^{(3)}_K\iy^{(3)})\vv\\
    =\;& \sum_{\substack{j_1^*,\dots,j_d^* \in [n]\\ j_1,\dots,j_t \in [n],\; j'_1,\dots,j'_t \in [n]\\i_1,\dots,i_t \in [m],\;i'_1,\dots,i'_t \in [m]}}\frac{\ell_{\iy^{(1)}}\left( \prod_{s=1}^d \pi^{(1)}_{i^*_sj_s^*}\prod_{s=1}^t\pi^{(1)}_{i_sj_s}\prod_{s=1}^t\pi^{(1)}_{i'_sj'_s}\right)}{\prod_{s=1}^d\beta_{j^*_s}\prod_{s=1}^t\beta_{j_s}\prod_{s=1}^t\beta_{j'_s}}\\
    \cdot& \sum_{\substack{j''_1,\dots,j''_t \in [n]\\ k''_1,\dots,k''_t \in [p]}}g(\prod_{s=1}^t\pi^{(2)}_{j''_sk''_s},\prod_{s=1}^t\pi^{(1)}_{i_sj_s})g(\prod_{s=1}^t\pi^{(2)}_{j''_sk''_s},\prod_{s=1}^t\pi^{(1)}_{i'_sj'_s})\\
    =\;& \sum_{\substack{j_1^*,\dots,j_d^* \in [n]\\ j''_1,\dots,j''_t \in [n]\\ k''_1,\dots,k''_t \in [p]}}\frac{1}{\prod_{s=1}^d\beta_{j^*_s}}\sum_{\substack{j_1,\dots,j_t \in [n],\; j'_1,\dots,j'_t \in [n]\\i_1,\dots,i_t \in [m],\;i'_1,\dots,i'_t \in [m]}}\ell_{\iy^{(1)}}\left( \prod_{s=1}^d \pi^{(1)}_{i^*_sj_s^*}\prod_{s=1}^t\pi^{(1)}_{i_sj_s}\prod_{s=1}^t\pi^{(1)}_{i'_sj'_s}\right)\\
    \cdot & \frac{g(\prod_{s=1}^t\pi^{(2)}_{j''_sk''_s},\prod_{s=1}^t\pi^{(1)}_{i_sj_s})}{\prod_{s=1}^t\beta_{j_s}}\frac{g(\prod_{s=1}^t\pi^{(2)}_{j''_sk''_s},\prod_{s=1}^t\pi^{(1)}_{i'_sj'_s})}{\prod_{s=1}^t\beta_{j'_s}}
\end{align*}
Let $I= \{(i^*_1,j^*_1),\dots,(i^*_d,j^*_d)\}$.  Using the fact that $\overline{\M}_t(e_I^{(1)}\iy^{(1)}) \succeq 0$, we have 
    \begin{multline*}
        \sum_{\substack{j_1,\dots,j_t \in [n],\; j'_1,\dots,j'_t \in [n]\\i_1,\dots,i_t \in [m],\;i'_1,\dots,i'_t \in [m]}}\ell_{\iy^{(1)}}\left( \prod_{s=1}^d \pi^{(1)}_{i^*_sj_s^*}\prod_{s=1}^t\pi^{(1)}_{i_sj_s}\prod_{s=1}^t\pi^{(1)}_{i'_sj'_s}\right)\\
    \cdot  \frac{g(\prod_{s=1}^t\pi^{(2)}_{j''_sk''_s},\prod_{s=1}^t\pi^{(1)}_{i_sj_s})}{\prod_{s=1}^t\beta_{j_s}}\frac{g(\prod_{s=1}^t\pi^{(2)}_{j''_sk''_s},\prod_{s=1}^t\pi^{(1)}_{i'_sj'_s})}{\prod_{s=1}^t\beta_{j'_s}} \geq 0.
    \end{multline*}
    This implies that $\overline{\M}_t(e_K^{(3)}\iy^{(3)}) \succeq 0$ for all $K \subset [p \times m]$, which completes the proof.
\end{proof}

\begin{proposition}\label{lemma triangular inequality}
Let $\XX$, $\YY$ and $\ZZ$ be metric measure spaces.  We then have the following
\begin{equation*}
\Delta_{p,r}(\XX,\ZZ) \leq \Delta_{p,r}(\XX,\YY) + \Delta_{p,r}(\YY,\ZZ).
\end{equation*}
\end{proposition}

\begin{proof}[Proof of Proposition \ref{lemma triangular inequality}]

Let $\iy^{(1)} \in \Pi_r(\alpha,\beta)$ and $\iy^{(2)} \in \Pi_r(\beta,\gamma)$ be minimizers of $\eqref{hierarchy Schmudgen discrete GW}(\XX,\YY)$ and $\eqref{hierarchy Schmudgen discrete GW}(\YY,\ZZ)$, respectively -- these exist because the sets $\Pi_r(\alpha,\beta)$ and $\Pi_r(\beta,\gamma)$ are compact.  Define the measures $\mu \in (\CX \times \CY)^2$ and $\nu \in (\CY \times \CZ)^2$ by
\begin{equation*}
\mu(x_i,y_j,x_{i'},y_{j'}) = \ell_{\iy^{(1)}}(\pi^{(1)}_{ij}\pi^{(1)}_{i'j'}),\quad \text{and} \quad \nu(y_j,z_{k},y_{j'},z_{k'}) = \ell_{\iy^{(2)}}(\pi^{(2)}_{jk}\pi^{(2)}_{j'k'}). 
\end{equation*}
We then have
\begin{equation*}
\Delta_{p,r}(\XX,\YY)= \Big (\int_{(\CX \times \CY)^2}|d_{\CX}(x,x')-d_{\CY}(y,y')|^pd\mu \Big)^{1/p}, 
\end{equation*}
and 
\begin{equation*}
\Delta_{p,r}(\YY,\ZZ)= \Big (\int_{(\CY \times \CZ)^2}|d_{\CY}(y,y')-d_{\CZ}(z,z')|^pd\nu \Big)^{1/p}.
\end{equation*}

Construct $\iy \in \RR^{s(mnp,2r)}$ based on the conclusions on the gluing lemma (Lemma \ref{thm:gluinglemma}).  First, it is straightforward to see that $\eta$, defined by 
\begin{equation*}
\eta(x_i,y_j,z_k,x_{i'},y_{j'},z_{k'}) = \ell_{\iy}(\pi_{ijk}\pi_{i'j'k'}),
\end{equation*}
specifies a probability measure on $(\CX \times \CY \times \CZ)^2$.  Second, based on the conclusions of the gluing lemma, the marginal measure of $\eta$ on $(\CX \times \CY)^2$ and $(\CY \times \CZ)^2$ are $\mu$ and $\nu$, respectively.  In addition, consider $\iy^{(3)}$, which also specifies another probability measure $\xi$ on $(\CZ \times \CX)^2$, and is the marginal measure of $\eta$ on $(\CZ \times \CX)^2$.  We then have
\begin{equation*}
\begin{aligned}
& \Delta_{p,r}(\XX,\ZZ)\\
\leq \;& \Big (\int_{(\CZ \times \CX)^2}|d_{\CZ}(z,z')-d_{\CX}(x,x')|^pd\xi \Big)^{1/p}\\
=\;& \Big(\int_{(\CX \times \CY \times \CZ)^2}|d_{\CZ}(z,z')-d_{\CY}(y,y')+d_{\CY}(y,y')-d_{\CX}(x,x')|^pd\eta \Big)^{1/p}\\
\leq\; & \Big(\int_{(\CX \times \CY \times \CZ)^2}|d_{\CZ}(z,z')-d_{\CY}(y,y')|^pd\eta \Big)^{1/p}+\Big(\int_{(\CX \times \CY \times \CZ)^2}|d_{\CY}(y,y')-d_{\CX}(x,x')|^pd\eta \Big)^{1/p}\\
=\;& \Big(\int_{(\CY \times \CZ)^2}|d_{\CY}(y,y')-d_{\CZ}(z,z')|^pd\nu \Big)^{1/p}+ \Big(\int_{(\CX \times \CY)^2}|d_{\CX}(x,x')-d_{\CY}(y,y')|^pd\mu \Big)^{1/p}\\
=\;& \Delta_{p,r}(\YY,\ZZ)+\Delta_{p,r}(\XX,\YY).
\end{aligned}
\end{equation*}
\end{proof}


\section{Numerical experiments}

In this section, we apply the first level of the moment-SOS hierarchy to solve instances of the Gromov-Wasserstein problem.  We perform these experiments on data generated from a computer vision dataset.  Our experiments serve a number of purposes.  First, we show that the first level of the hierarchy already provides a strong relaxation -- in many instances, we are able to compute {\em globally optimal} solutions, together with a {\em proof} of global optimality.  Second, we are able to solve moderate-sized instances in a reasonable amount of time.  In the largest problem instance, we were able to solve a GW instance with $50$ atoms in the source as well as the target to {\em global optimality}, all within a one-hour time limit on a laptop.  The corresponding SDP is of size $50^2$ by $50^2$.  

In what follows, we focus on computing solutions to \eqref{eq:lp_definition} in the case where $p=2$, and where the metric $d$ is the Euclidean distance.  The source and target are discrete metric measure spaces $(\CX,d_{\CX},\mu_{\CX})$, with dimension $m$, and $(\CY,d_{\CY},\mu_{\CY})$, with $n$.  The first order Schmüdgen-type moment relaxation is given as the following SDP form:
 \begin{equation} \label{numerical Schmudgen discrete GW} 
\begin{aligned}
\Delta_{2,1}(\CX,\CY):=\min \qquad & \sum_{i,k \in [m],\; j,l \in [n]}L_{ij,kl}\cdot \iy_{ij,kl} \\
\mbox{subject to} \qquad & \iy \in \RR^{s(mn,2)},\;\iy_0 =1, \;\\ 
&\M_1(\iy) \succeq 0,\; \iy_{ij,kl} \geq 0 \; \forall i, k \in [m],\ j,l \in [n],\\
&\M_{0}(r_i(\pi) \iy) = 0, \;\M_{0}(r_{i_1}e_{i_2j}(\pi) \iy) =0\; \forall i, i_1, i_2 \in [m],\ j \in [n]\\
&\M_{0}(c_j(\pi) \iy) = 0, \;\M_{0}(c_{j_1}e_{ij_2}(\pi) \iy) =0\; \forall i\in [m],\;  j, j_1,j_2 \in [n].
\end{aligned} \tag{S-DGW-1}
\end{equation}
We solve the associated SDPs using the specialized software RNNAL \cite{hou2025lowrankaugmentedlagrangianmethod}.  The software is designed specifically for SDPs whereby the matrix variable is required to be both PSD and non-negative, also known as a {\em doubly non-negative} (DNN) matrix.  The SDP instances \eqref{numerical Schmudgen discrete GW} is in fact a DNN relaxation.  The RNNAL software is built using the ideas of the low-rank augmented Lagrangian method \cite{hou2025lowrankaugmentedlagrangianmethod}.


\textbf{Verifying Exactness.}  The standard method within the POP literature to determine exactness of an SDP relaxation is via a concept known as {\em flat extension}, in which one computes the rank of the moment matrices at the optimal solution of the primal SDP relaxation -- see, for instance \cite{lasserre2009moments}.  However, in practice, the rank condition is difficult to verify because of numerical error.

In this paper, we use two alternative quantities to indicate exactness.  The first quantity is the ratio between the largest eigenvalue and the second largest eigenvalue of the moment matrix of the optimal solution -- this is an approach that is widely used; see, for instance  \cite{burer2024slightly} and references therein.  
More specifically, we compute the following {\em eigenvalue ratio}:
\begin{equation*}
    \mbox{eigenvalue ratio} := \frac{\lambda_2(\M_1(\iy^*))}{\lambda_1(\M_1(\iy^*))},
\end{equation*}
where $\lambda_1(\M_1(\iy^*))$ and $\lambda_2(\M_1(\iy^*))$ denote the largest eigenvalue and the second largest eigenvalue of $\M_1(\iy^*)$, and $\iy^*$ denotes the optimal solution of~\ref{numerical Schmudgen discrete GW}.  If the ratio is close to zero, we know that the matrix $\M_1(\iy^*)$ is close to being rank-one, which indicates that the relaxation is exact.

The second quantity we compute measures how close the lower bound (defined as the optimal value of~\eqref{numerical Schmudgen discrete GW}) to the GW distance, denoted by {\em error ratio}, defined as:
\begin{equation*}
    \mbox{error ratio} := \frac{\sum_{i,k \in [m],\; j,l \in [n]}L_{ij,kl}\cdot \iy^*_{ij}\iy^*_{kl}}{\Delta_{2,1}(\CX,\CY)}.
\end{equation*}
Because of the marginal condition~\eqref{marginal conditions}, the coupling $(\iy^*_{ij})_{i \in [m], j\in [n]}$ defines a feasible coupling in $\Pi(\mu_{\CX},\mu_{\CY})$ -- here $\iy^*_{ij}$ is the component of $\iy^*$ indexed by the monomial $\pi_{ij}$. Thus, $\sum_{i,k \in [m],\; j,l \in [n]}L_{ij,kl}\cdot \iy^*_{ij}\iy^*_{kl}$ defines an upper bound on the GW distance while $\Delta_{2,1}(\CX,\CY)$ is a lower bound of the GW distance. As a result, the {\em error ratio} is at least one; in particular, if the ratio equals to one, it indicates that the relaxation \eqref{numerical Schmudgen discrete GW} is exact.  As such, we will be measuring how close the ratio is to one.

\textbf{Implementation.}  All the experiments are performed on a Macbook with the Apple M1 Pro chip with 8 cores and 16GB RAM.  We distinguish between two sets of experiments:  In the first instance, we solve GW instances where the source and target spaces are atomic measures with equal number of atoms.  In the second instance, the number of atoms in the source and the target distribution are different.  The size of PSD relaxation -- specifically, the matrix dimension -- is given by {\tt N = m$\times$n}.  In the following, the parameters {\tt lower bound, time, eig. ratio}, and \texttt{err. ratio} refer to the optimal value of~\eqref{numerical Schmudgen discrete GW}, the compute time (in seconds), the eigenvalue ratio, and the error ratio respectively. The tolerance in SDP solver is set to $10^{-6}$.  We declare that a GW instance is {\em solved} if the error ratio is at most $1.0001$ and if the eigenvalue ratio is smaller than $10^{-4}$.  Our criteria is similar to those used in~\cite{burer2024slightly,consolini2023sharp,eltved2023strengthened}. 

In the first set of experiments, we sample the same number of points $m$ on the surface of a cat across three different poses.  We start with $m=5$, and we increase the number of points in multiples of $5$ until we reach $50$.  We set a run-time limit of one-hour (the run-time exceeded the one-hour limit for problem instances beyond $m=50$ points).  For each $m$ and each pair of poses, we compute the GW distance between the corresponding set of $m$ points, and we present the full set of results in Table~\ref{tab: cat}.
\begin{figure}[H]
\centering
\begin{subfigure}[b]{0.3\textwidth}
    \centering
    \includegraphics[width = 0.5\textwidth]{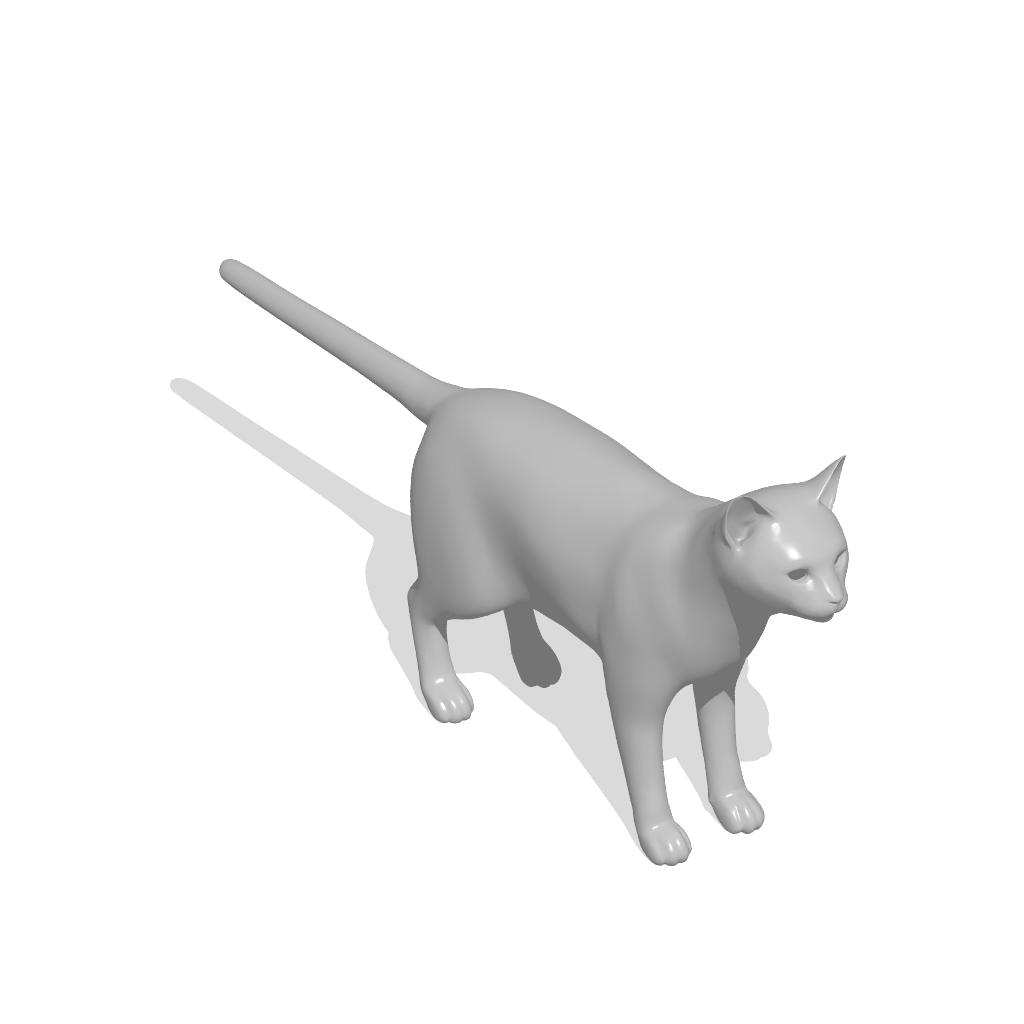}
\end{subfigure}
\begin{subfigure}[b]{0.3\textwidth}
    \centering
    \includegraphics[width = 0.5\textwidth]{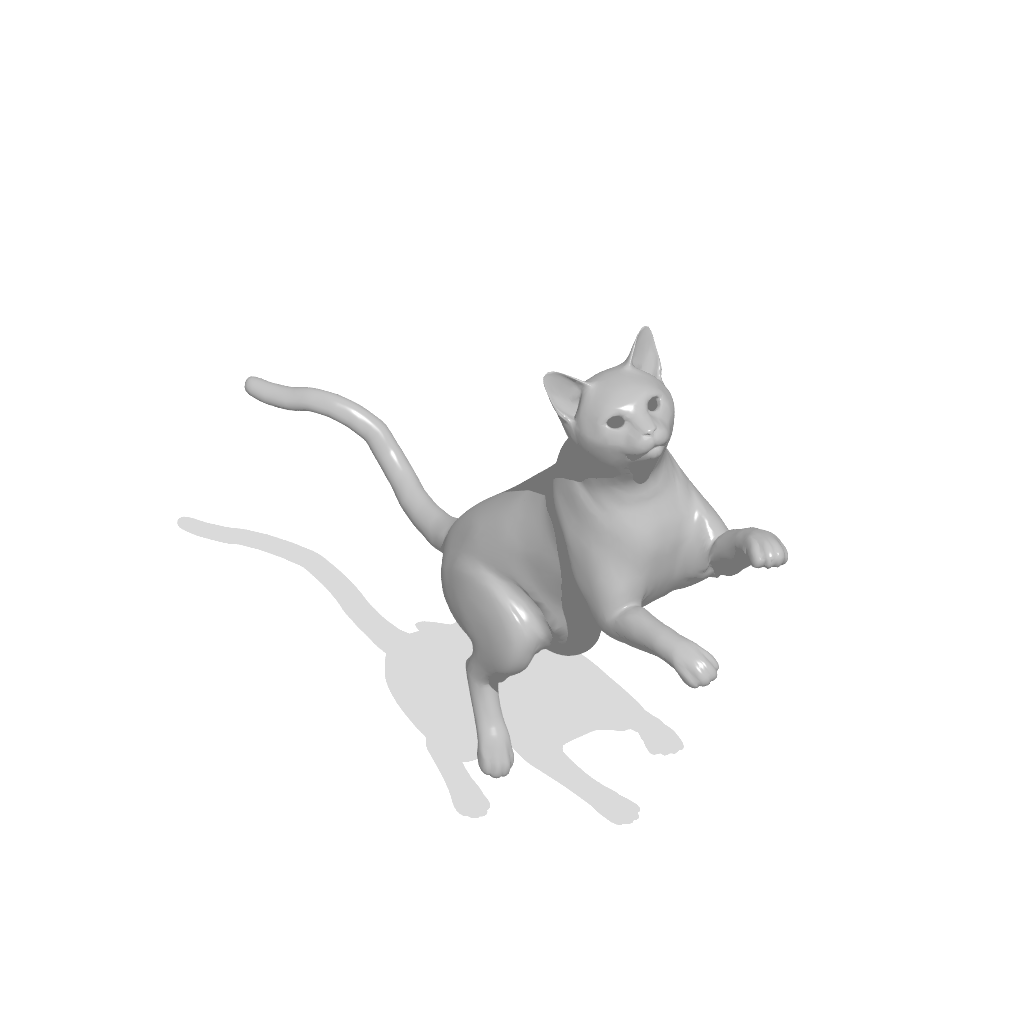}
\end{subfigure}
\begin{subfigure}[b]{0.3\textwidth}
    \centering
    \includegraphics[width = 0.5\textwidth]{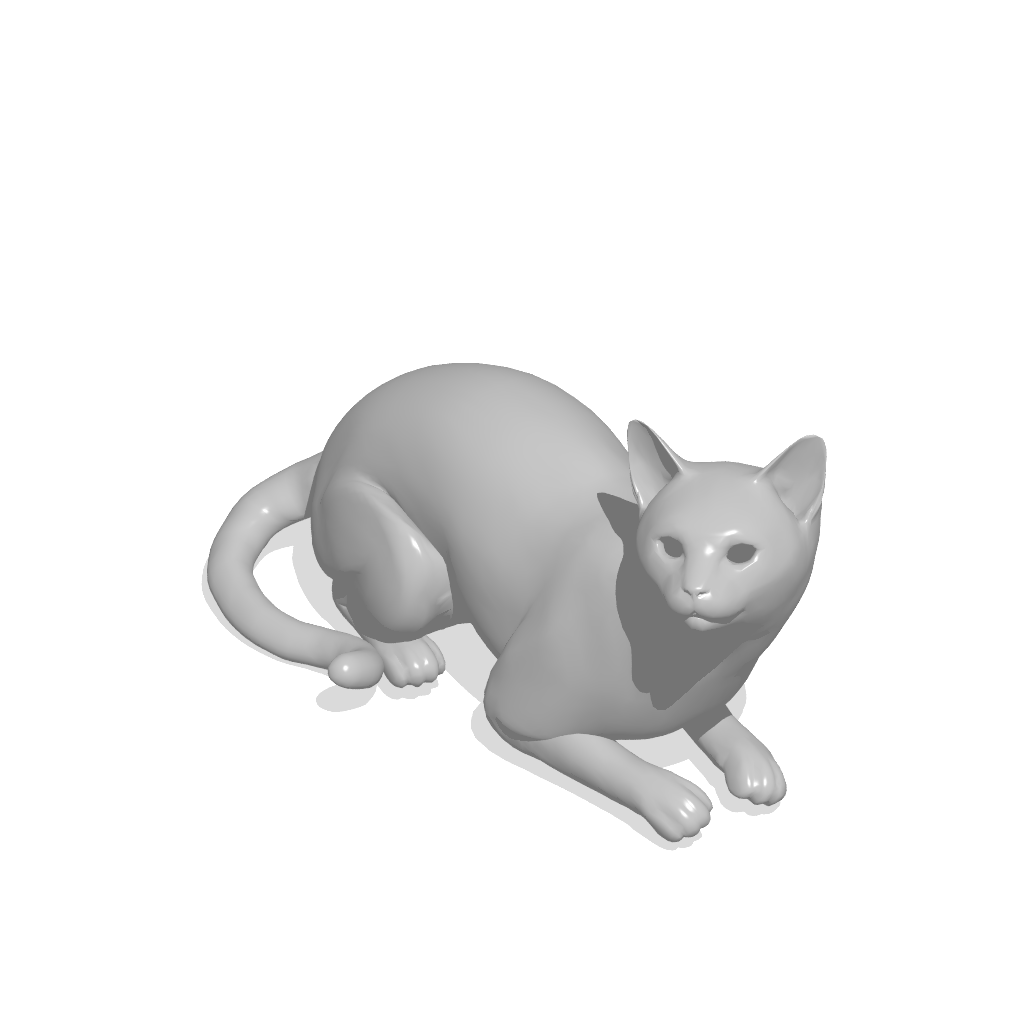}
\end{subfigure}
\caption{Selected poses of cats from the TOSCA dataset \cite{bronstein2008numerical}.}
\label{fig: cats}
\end{figure}

{\color{blue}{
\begin{tiny}
\begin{longtable}[c]{cccccc}
\caption{Computational results for DGW shape correspondence problems}\\
\toprule
cat problem & type & lower bound $(GW_1)$ & time & eig. ratio & err. ratio\\
\midrule
\endfirsthead

\multicolumn{6}{c}%
{{ Table \thetable\ continued from previous page}} \\
\toprule
cat problem & type & lower bound & time & eig. ratio & err. ratio\\
\endhead
\midrule
\multicolumn{6}{r}{{Continued on next page}} \\
\midrule
\endfoot
\bottomrule
\endlastfoot

\multirow{3}{*}{$N=25$}& $(1) \leftrightarrow (2)$ &4.18663e+02	&1.58975e-01 &1.04843e-08 &1.00000e+00 \\
                       & $(2) \leftrightarrow (3)$ &8.42575e+01 &3.24520e-02 &5.72044e-08 &1.00000e+00 \\
                       & $(3) \leftrightarrow (1)$ &3.36854e+02 &1.61330e-02 &3.80274e-08 &1.00000e+00\\[4pt] 
                       
\multirow{3}{*}{$N=100$}& $(1) \leftrightarrow (2)$ &1.24246e+02 &9.45294e-01 &1.14362e-10 &1.00000e+00 \\
                       & $(2) \leftrightarrow (3)$ &9.32959e+01 &2.73890e-01 &7.74060e-09 &1.00000e+00\\
                       & $(3) \leftrightarrow (1)$ &1.24422e+02 &1.36546e+00 &8.05832e-03 &1.00211e+00\\ [4pt] 

\multirow{3}{*}{$N=225$}& $(1) \leftrightarrow (2)$ &1.15611e+02 &6.56165e+00 &7.45225e-05 &1.00004e+00\\
                       & $(2) \leftrightarrow (3)$ &1.90305e+02 &8.10347e+00 &3.41520e-04 &1.00337e+00 \\
                       & $(3) \leftrightarrow (1)$ &1.92000e+02 &3.72385e+00 &9.23883e-04 &1.00058e+00\\ [4pt] 

\multirow{3}{*}{$N=400$}& $(1) \leftrightarrow (2)$ &9.86495e+01 &2.53831e+01 &1.29279e-05 &1.00001e+00\\
                       & $(2) \leftrightarrow (3)$ &1.08252e+02 &1.37160e+01 &7.14956e-03 &1.00097e+00\\
                       & $(3) \leftrightarrow (1)$ &1.01634e+02 &2.04483e+01 &1.09280e-04 &1.00012e+00\\ [4pt]

\multirow{3}{*}{$N=625$}& $(1) \leftrightarrow (2)$ &1.34882e+02 &2.79877e+02 &2.40831e-04 &1.00023e+00\\
                       & $(2) \leftrightarrow (3)$ &9.17491e+01 &6.67915e+01 &2.93618e-04 &1.00001e+00\\
                       & $(3) \leftrightarrow (1)$ &8.02853e+01 &5.04674e+01 &6.15517e-03 &1.00025e+00\\ [4pt]

\multirow{3}{*}{$N=900$}& $(1) \leftrightarrow (2)$ &9.25423e+01 &3.52483e+02 &2.92345e-03 &1.00056e+00\\
                       & $(2) \leftrightarrow (3)$ &1.10631e+02  &1.89159e+02 &5.00012e-03 &1.00010e+00\\
                       & $(3) \leftrightarrow (1)$ &8.21558e+01 &2.07985e+02 &1.70306e-03 &1.00012e+00\\ [4pt] 

\multirow{3}{*}{$N=1225$}& $(1) \leftrightarrow (2)$ &7.71375e+01 &5.28732e+02 &1.27465e-02 &1.00979e+00\\
                       & $(2) \leftrightarrow (3)$ &9.92526e+01 &3.39351e+02  &2.85686e-02 &1.00029e+00\\
                       & $(3) \leftrightarrow (1)$ &8.33600e+01 &4.80225e+02 &1.55544e-04 &1.00027e+00\\ [4pt]

\multirow{3}{*}{$N=1600$}& $(1) \leftrightarrow (2)$ &6.35761e+01 &9.96248e+02 &1.67914e-03 &1.00017e+00\\
                       & $(2) \leftrightarrow (3)$ &8.71064e+01 &7.68950e+02 &1.04060e-03 &1.00002e+00\\
                       & $(3) \leftrightarrow (1)$ &6.62119e+01 &8.67505e+02 &8.80503e-05 &1.00000e+00\\ [4pt]
                       
\multirow{3}{*}{$N=2025$}& $(1) \leftrightarrow (2)$ &7.94103e+01 &2.60458e+03 &2.42851e-02 &1.00128e+00\\
                       & $(2) \leftrightarrow (3)$ &8.87060e+01  &1.21919e+03 &2.73319e-06 &1.00000e+00\\
                       & $(3) \leftrightarrow (1)$ &8.15298e+01 &1.36268e+03 &1.20835e-02 &1.00011e+00 \\ [4pt]

\multirow{3}{*}{$N=2500$}& $(1) \leftrightarrow (2)$ &9.66241e+01 &2.39597e+03 &2.05495e-02 &1.00043e+00 \\
                       & $(2) \leftrightarrow (3)$ &6.86050e+01 &2.10420e+03 &8.09876e-04 &1.00009e+00   \\ 
                       & $(3) \leftrightarrow (1)$ &9.16975e+01 &2.92626e+03 &2.62067e-03 &1.00304e+00\\ [4pt]  
\label{tab: cat}
\end{longtable}
\end{tiny}
}}

\begin{figure}[H]
\centering
\begin{subfigure}[b]{0.3\textwidth}
    \centering
    \includegraphics[width = 0.5\textwidth]{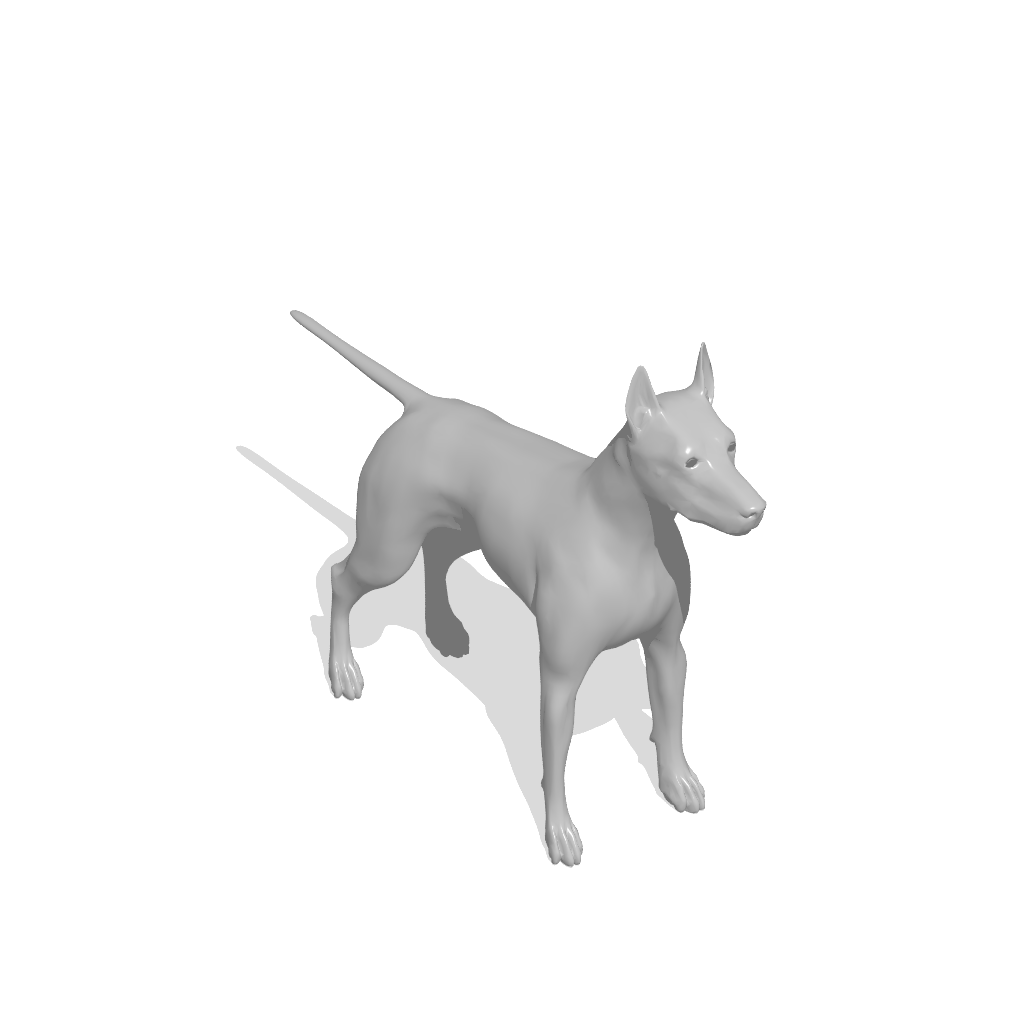}
\end{subfigure}
\begin{subfigure}[b]{0.3\textwidth}
    \centering
    \includegraphics[width = 0.5\textwidth]{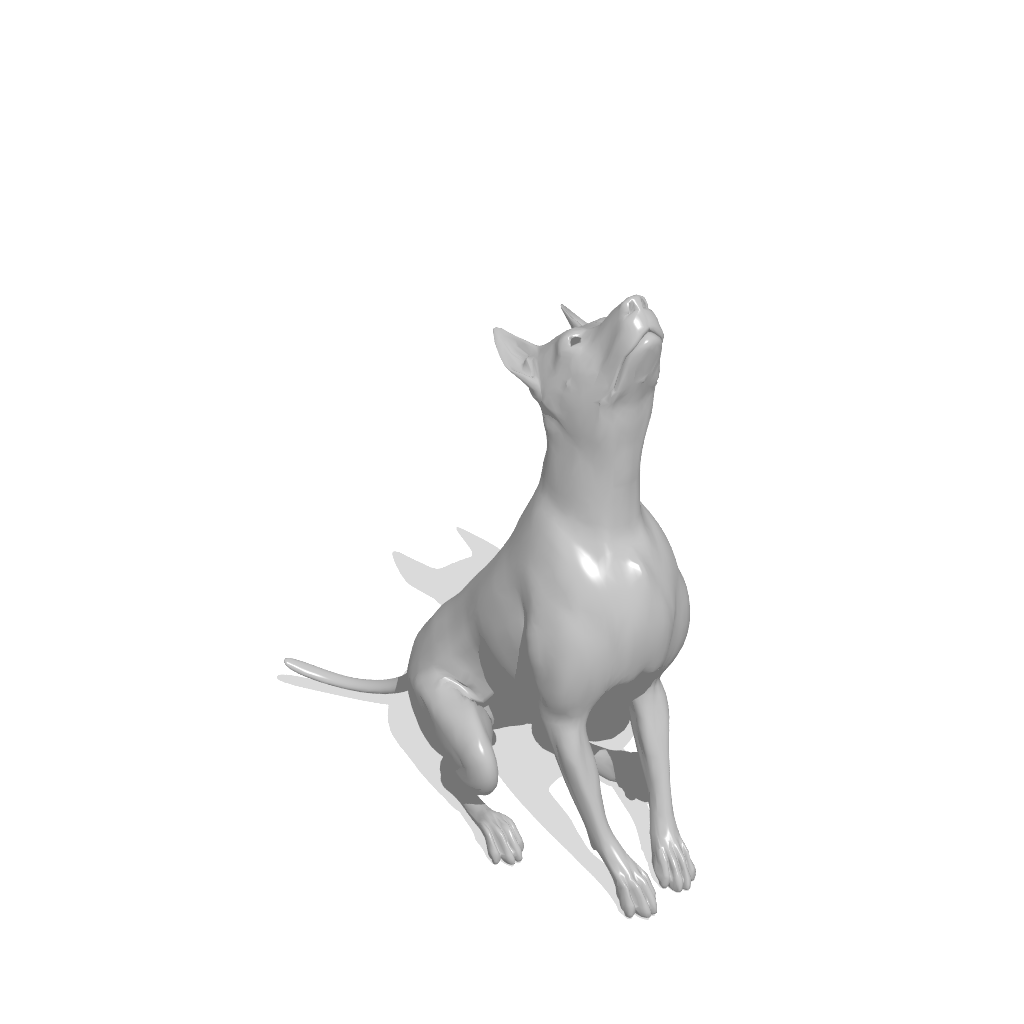}
\end{subfigure}
\begin{subfigure}[b]{0.3\textwidth}
    \centering
    \includegraphics[width = 0.5\textwidth]{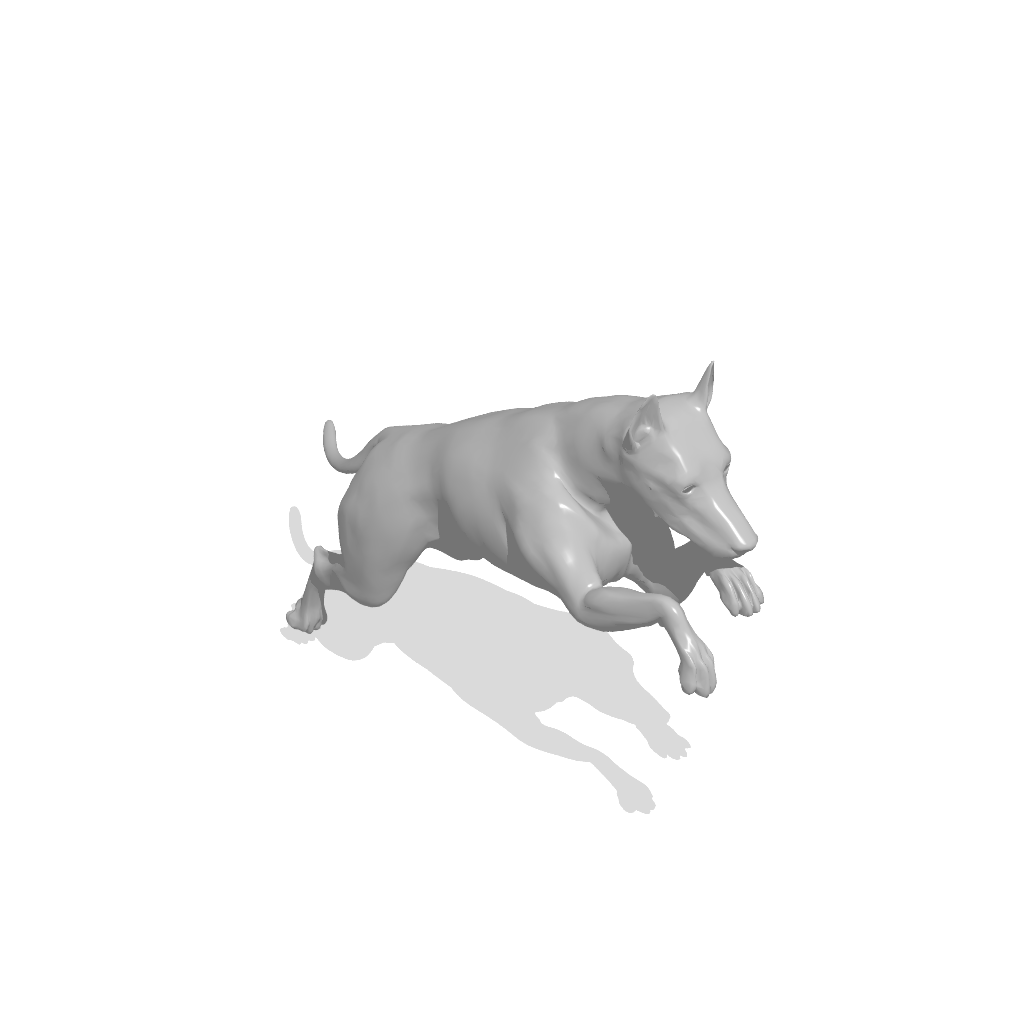}
\end{subfigure}
\caption{Selected poses of dogs from the TOSCA dataset \cite{bronstein2008numerical}.}
\label{fig: dogs}
\end{figure}

\begin{figure}[H]
\centering
\begin{subfigure}[b]{0.3\textwidth}
    \centering
    \includegraphics[width = 1.0\textwidth]{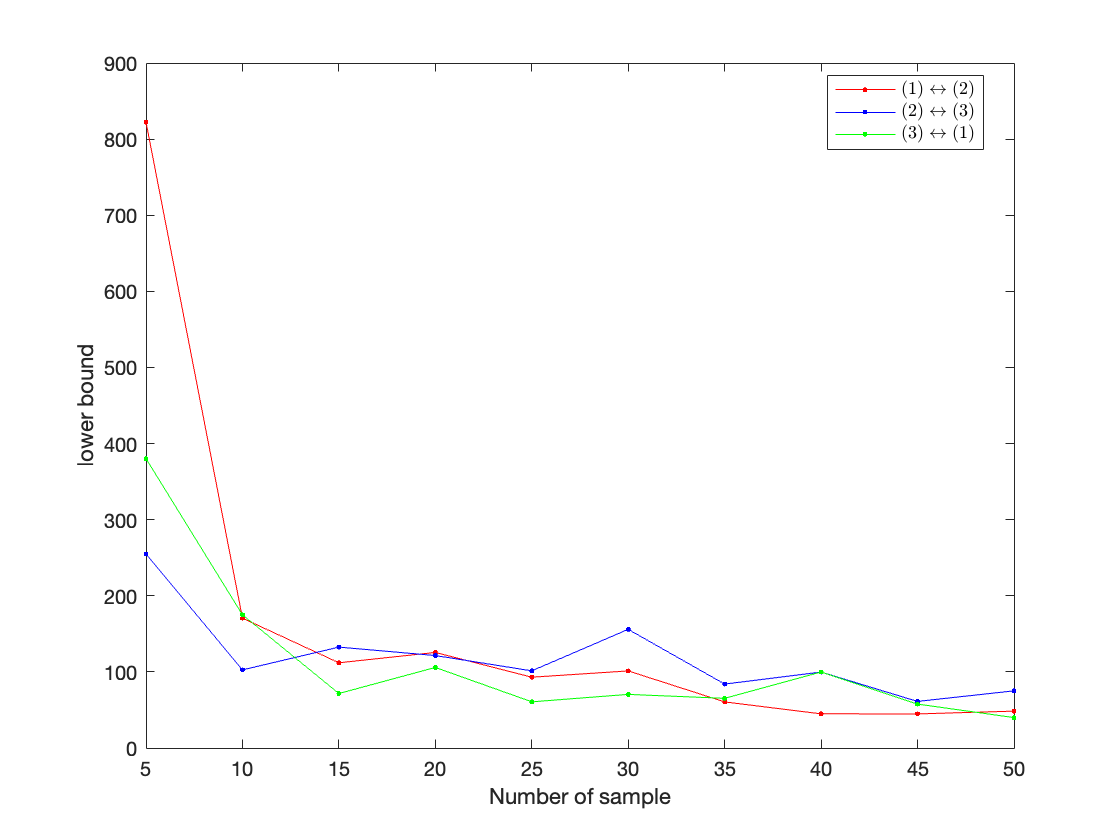}
\end{subfigure}
\begin{subfigure}[b]{0.3\textwidth}
    \centering
    \includegraphics[width = 1.0\textwidth]{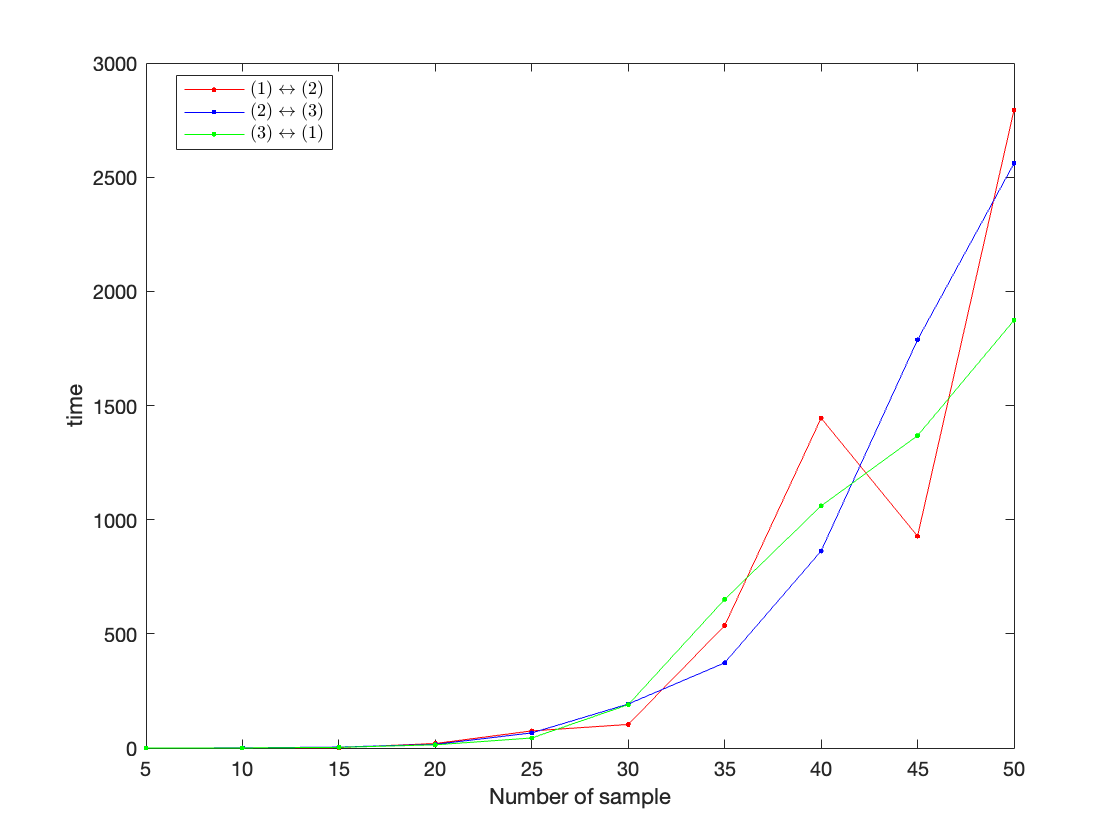}
\end{subfigure}
\begin{subfigure}[b]{0.3\textwidth}
    \centering
    \includegraphics[width = 1.0\textwidth]{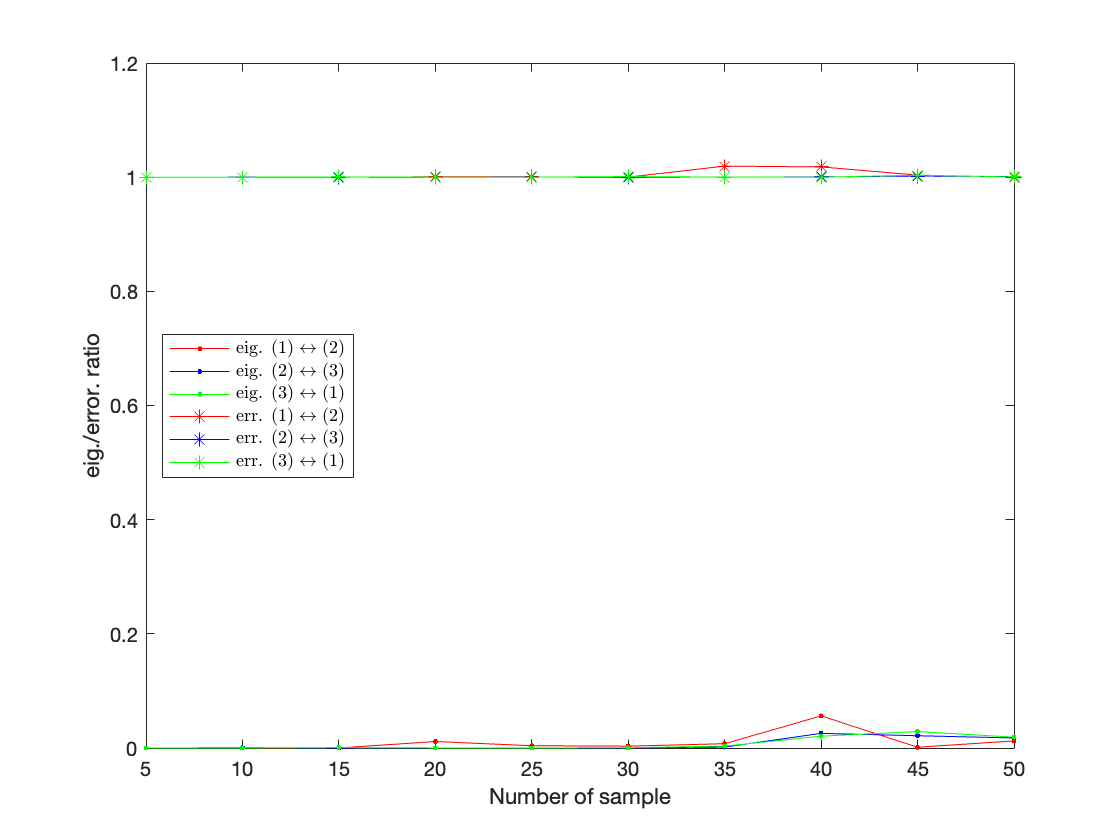}
\end{subfigure}
\caption{Results of the dog-shapes experiments}
\label{results: dog}
\end{figure}

\begin{figure}[H]
\centering
\begin{subfigure}[b]{0.3\textwidth}
    \centering
    \includegraphics[width = 0.5\textwidth]{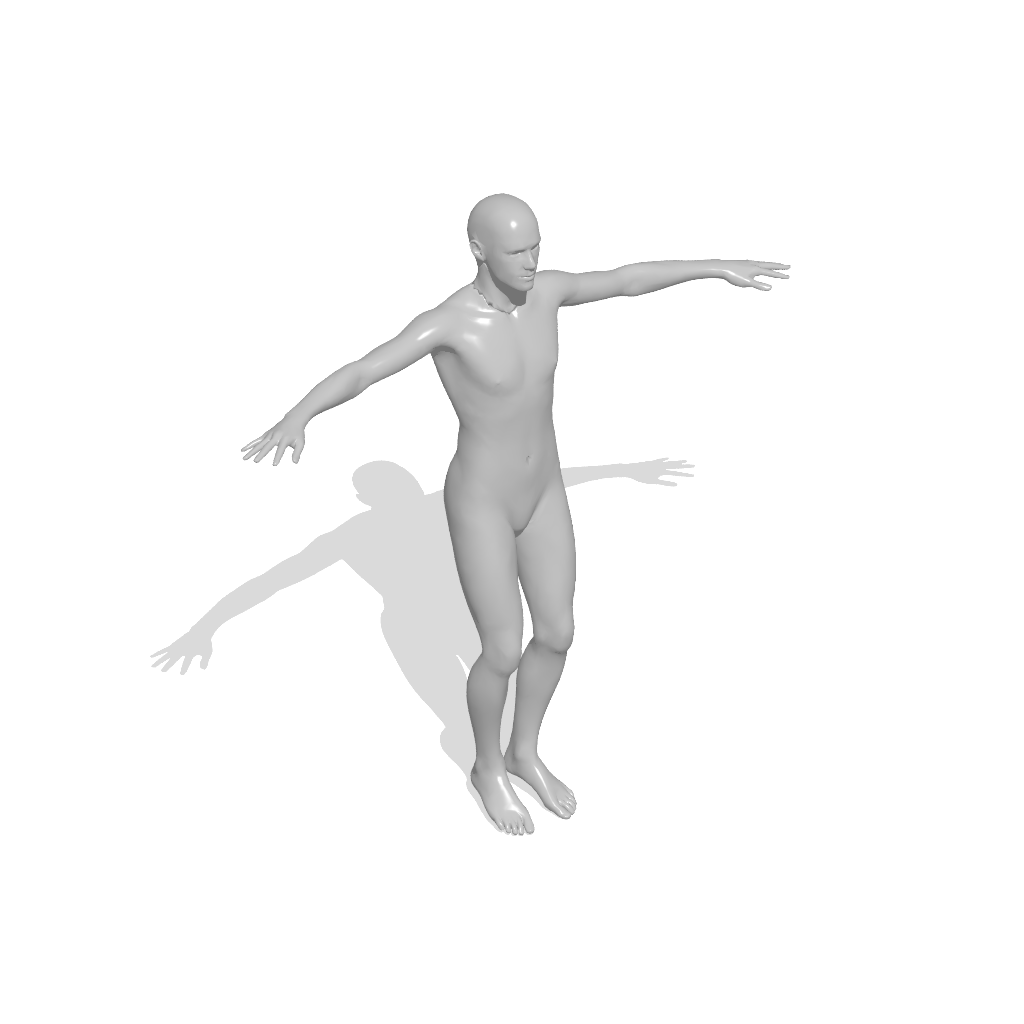}
\end{subfigure}
\begin{subfigure}[b]{0.3\textwidth}
    \centering
    \includegraphics[width = 0.5\textwidth]{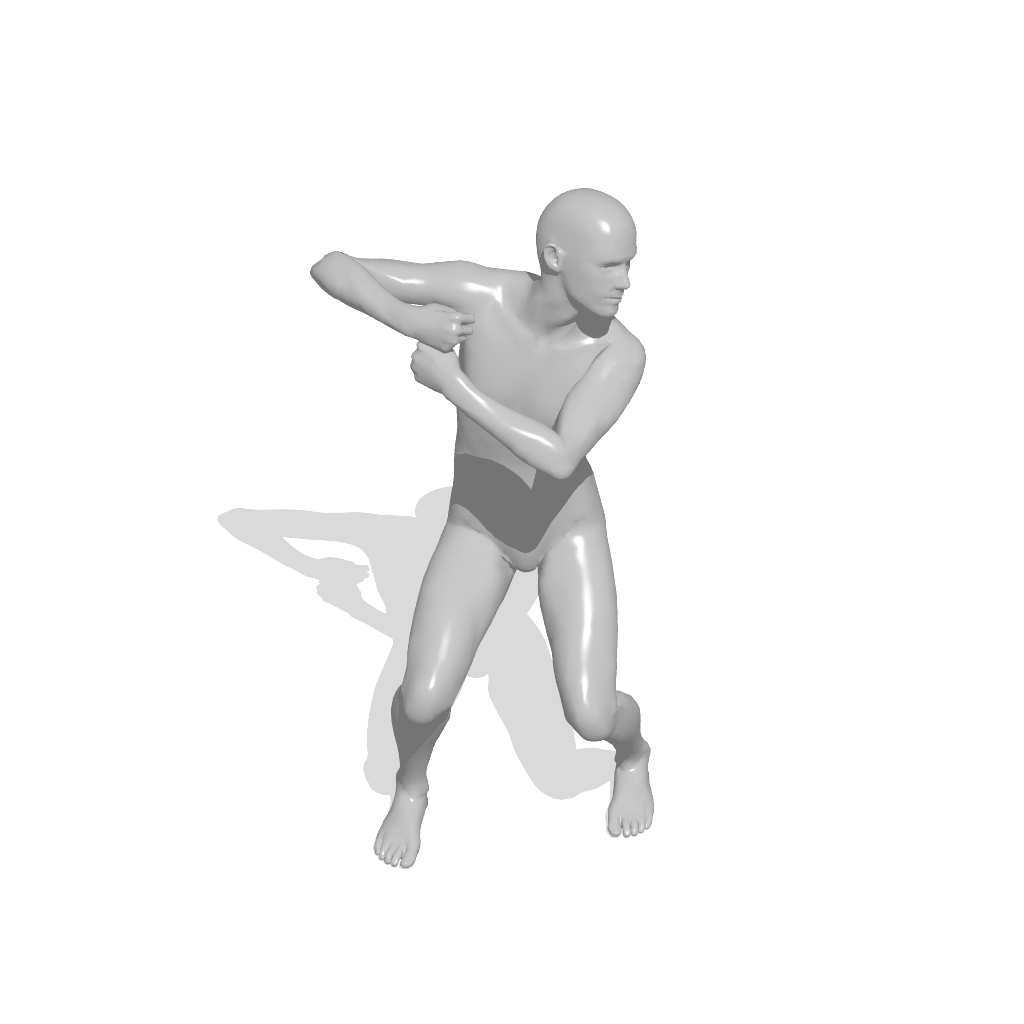}
\end{subfigure}
\begin{subfigure}[b]{0.3\textwidth}
    \centering
    \includegraphics[width = 0.5\textwidth]{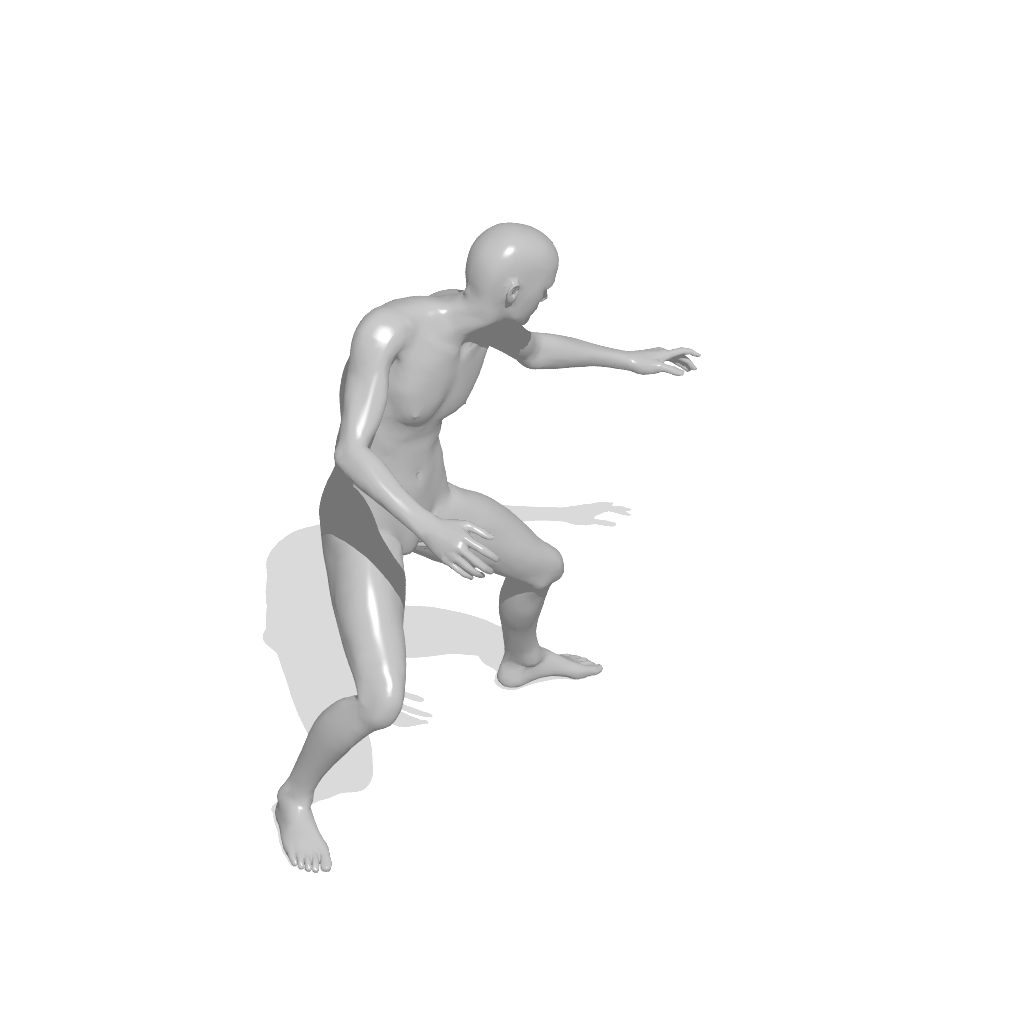}
\end{subfigure}
\caption{Selected poses of a human from the TOSCA dataset \cite{bronstein2008numerical}.}
\label{fig: humans}
\end{figure}

\begin{figure}[H]
\centering
\begin{subfigure}[b]{0.3\textwidth}
    \centering
    \includegraphics[width = 1.0\textwidth]{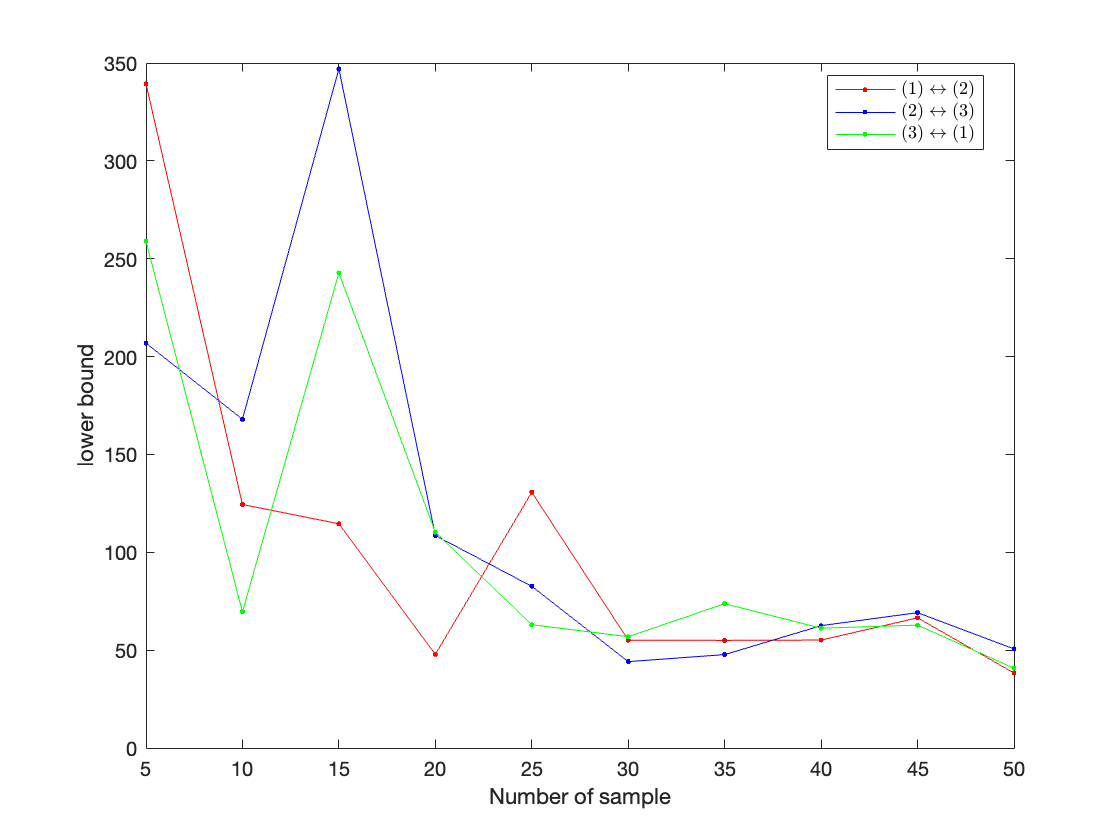}
\end{subfigure}
\begin{subfigure}[b]{0.3\textwidth}
    \centering
    \includegraphics[width = 1.0\textwidth]{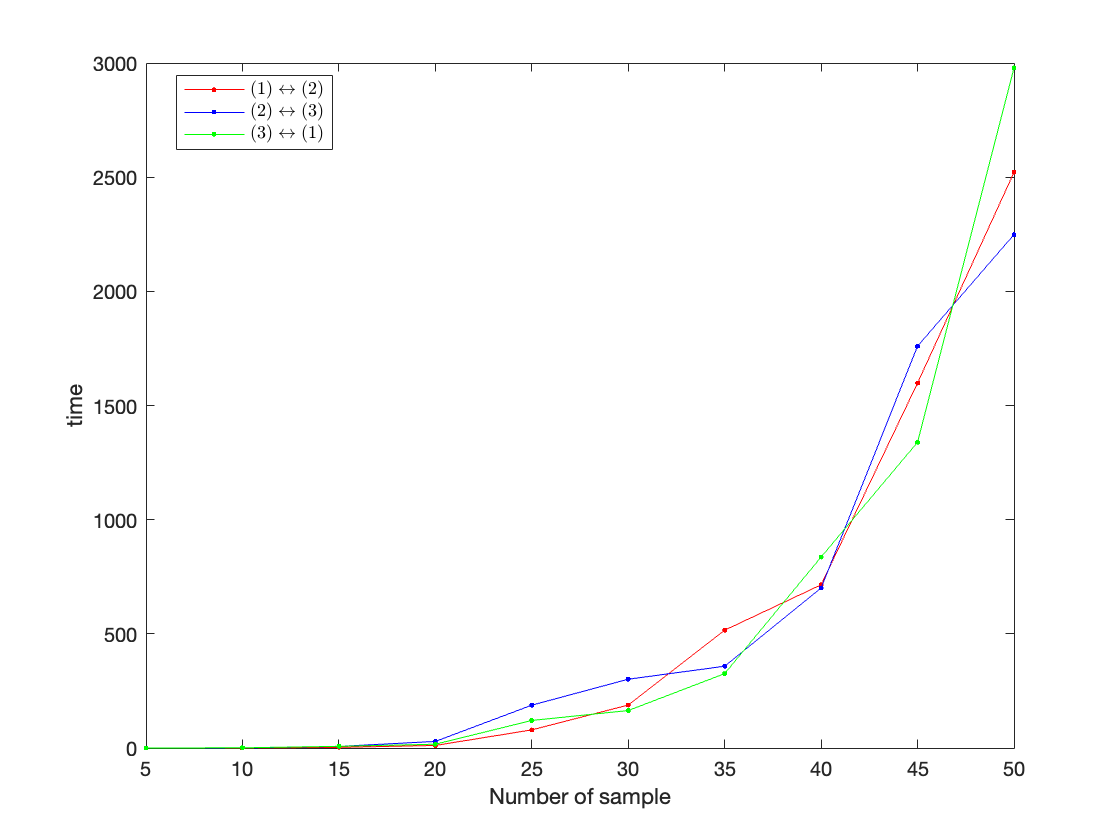}
\end{subfigure}
\begin{subfigure}[b]{0.3\textwidth}
    \centering
    \includegraphics[width = 1.0\textwidth]{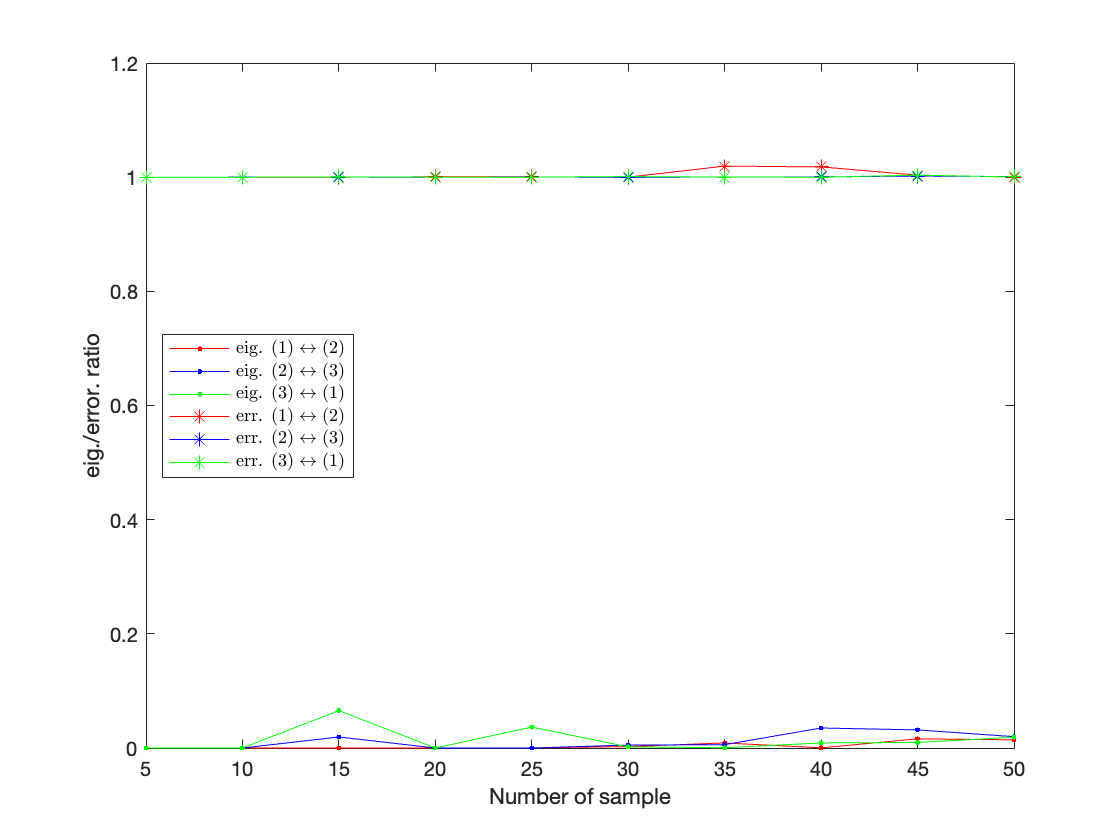}
\end{subfigure}
\caption{Results of the human-shapes experiments}
\label{results: human}
\end{figure}

Our results in Table~\ref{tab: cat} shows that the first level of the Schmüdgen-type moment relaxation is sufficient to solve the GW instances to global optimality for most instances, and within the one-hour time limit.  As a note, in the largest instance where $N=50$, the associated SDP has matrix dimension $50^2 \times 50^2$.  We repeat the same experimental set-up for the dog and human object --  see Figures \ref{fig: dogs} and \ref{fig: humans} -- and we present a summarized version of these results in Figures~\ref{results: dog} and~\ref{results: human}.  Our results are similar -- we observe exactness with similar frequencies these other objects.  

In the second set of experiments, we allow the number of atoms in the source and target distribution to be {\em different}.  For each object, we draw a different number of sample points for each pose.  The results of this experimental set-up is shown in Table \ref{tab: results for diff number of points}.

We observe a number of qualitative differences.  First, we notice that exactness holds in many instances, but the frequency at which exactness holds appears fewer than if the number of atoms in the source and target were equal.  Second, we notice that the compute time is substantially longer.  In some instances, particularly when there is a larger mismatch between the number of atoms, the numerical solver ran out of time.  We believe that there is a fundamental difference in the geometric structure of these relaxations between the two cases we consider that explains the difference in the observed strength of these relaxations.  We think investigating these questions can form the basis of interesting follow-up work.


{\color{blue}{
\begin{tiny}
\begin{longtable}[c]{c|cccccc}
\caption{Computational results for DGW shape correspondence problems}\\
\toprule
\multicolumn{2}{c}{problem} & type & lower bound $(GW_1)$ & time & eig. ratio & err. ratio\\
\midrule
\endfirsthead

\multicolumn{7}{c}%
{{ Table \thetable\ continued from previous page}} \\
\toprule
\multicolumn{2}{c}{problem} & type & lower bound & time & eig. ratio & err. ratio\\
\endhead
\midrule
\multicolumn{7}{r}{{Continued on next page}} \\
\midrule
\endfoot
\bottomrule
\endlastfoot

\multirow{9}{*}{\makecell{$m=10$,\\$n=15$,\\ $p=20$}} &\multirow{3}{*}{cat}
&$(1) \leftrightarrow (2)$   &3.15160e+02  &3.05801e+00	&8.64080e-05	&1.00000e+00 \\
&& $(2) \leftrightarrow (3)$ &1.66863e+02  &3.92271e+02	&1.96060e-04	&1.00020e+00 \\
&& $(3) \leftrightarrow (1)$ &- & \mbox{reach time limit} &- &-\\ [4pt]

&\multirow{3}{*}{dog}  
& $(1) \leftrightarrow (2)$  &3.52512e+02  &6.51612e+00	&4.63580e-05  &1.00000e+00\\
&& $(2) \leftrightarrow (3)$ &1.85722e+02  &1.63268e+02	&5.75250e-05  &1.00010e+00 \\
&& $(3) \leftrightarrow (1)$ &- & \mbox{reach time limit} &- &-\\ [4pt] 

&\multirow{3}{*}{human}
& $(1) \leftrightarrow (2)$ &2.91712e+02  &1.09873e+00	&1.11540e-07 &1.00000e+00 \\
&& $(2) \leftrightarrow (3)$ &1.40813e+02 &1.09959e+03	&1.72600e-01 &1.28610e+00 \\
&& $(3) \leftrightarrow (1)$ &- &\mbox{reach time limit} &- &-\\ [4pt] 

\multirow{9}{*}{\makecell{$m=20$,\\$n=25$,\\ $p=30$}} &\multirow{3}{*}{cat} 
&$(1) \leftrightarrow (2)$   &9.86789e+01 &3.27583e+02	&1.10370e-04 &1.00010e+00 \\
&& $(2) \leftrightarrow (3)$ &1.04749e+02 &3.03643e+03	&9.83000e-02 &1.06520e+00 \\
&& $(3) \leftrightarrow (1)$ &- & \mbox{reach time limit} &- &-\\ [4pt]

&\multirow{3}{*}{dog}  
& $(1) \leftrightarrow (2)$  &7.54419e+01	&1.20887e+02 &7.04650e-06	&1.00000e+00\\
&& $(2) \leftrightarrow (3)$ &9.63737e+01	&3.54485e+03 &9.80000e-03	&1.00680e+00 \\
&& $(3) \leftrightarrow (1)$ &- & \mbox{reach time limit} &- &-\\ [4pt] 

&\multirow{3}{*}{human}
& $(1) \leftrightarrow (2)$ &1.37548e+02  &4.90009e+02	&4.67760e-04	&1.00170e+00 \\
&& $(2) \leftrightarrow (3)$ &1.12621e+02 &6.24012e+02	&1.40000e-03	&1.00030e+00\\
&& $(3) \leftrightarrow (1)$ &- &\mbox{reach time limit} &- &-\\ [4pt] 

 \multirow{9}{*}{\makecell{$m=30$,\\$n=35$,\\ $p=40$}} &\multirow{3}{*}{cat}
 &$(1) \leftrightarrow (2)$   &7.15364e+01	&3.60362e+03 &1.72700e-01	&1.29670e+00 \\
&& $(2) \leftrightarrow (3)$ &6.33540e+01	&3.60611e+03  &1.99260e-04	 &1.00020e+00 \\
&& $(3) \leftrightarrow (1)$ &- & \mbox{reach time limit} &- &-\\ [4pt] 

&\multirow{3}{*}{dog}  
& $(1) \leftrightarrow (2)$  &6.72240e+01	&1.97285e+03 &1.06884e-03	&1.00165e+00\\
&& $(2) \leftrightarrow (3)$ &7.96791e+01	&3.60478e+03 &3.30028e-01	&1.46073e+00 \\
&& $(3) \leftrightarrow (1)$ &- & \mbox{reach time limit} &- &-\\ [4pt] 

&\multirow{3}{*}{human}
& $(1) \leftrightarrow (2)$ &7.28641e+01 &2.10029e+03 &2.30866e-03	&1.00433e+00 \\
&& $(2) \leftrightarrow (3)$ &6.30814e+01	&3.60153e+03	&4.90136e-02	&1.01479e+00 \\
&& $(3) \leftrightarrow (1)$ &- &\mbox{reach time limit} &- &-\\ [4pt]                                         
\label{tab: results for diff number of points}
\end{longtable}
\end{tiny}
}}
{\bf Key difference between our experimental set-up and the work in \cite{chen2023semidefinite,hou2025lowrankaugmentedlagrangianmethod}.}  It is necessary to point out that the work in \cite{hou2025lowrankaugmentedlagrangianmethod} also implements the RNNAL algorithm on a DNN-based relaxation of the GW problem -- the details can be found in Section 6.6 of \cite{hou2025lowrankaugmentedlagrangianmethod}.  In particular, the DNN relaxation in their work coincides with the first level of the Schmüdgen-type moment hierarchy.  We explain the differences between both sets of experiments.  First, our objectives are quite different.  In \cite{hou2025lowrankaugmentedlagrangianmethod}, the objective is to demonstrate numerical performance, specifically in the form of run-times.  In contrast, our objective is to show the strength of the relaxations.  As such, we report and focus on the eigenvalue ratio metrics as well as the error ratio metric as a measure of global optimality.  In particular, the work in \cite{hou2025lowrankaugmentedlagrangianmethod} does {\em not} discuss global optimality.  Second, the numerical meshes for the objects in the TOSCA dataset \cite{bronstein2008numerical} are generated independently.  From this perspective, the datasets used in both sets of experiments are technically different.

We also note that the numerical experiments in \cite{chen2023semidefinite} also explores exactness.  The scale of the experiments in \cite{chen2023semidefinite}, however, are far more modest.  We also investigate settings where the source and the target spaces have different number of atoms, which the work in \cite{chen2023semidefinite} does not.


\section{Conclusions and Future Directions}

In this paper, we developed tractable semidefinite relaxations for computing globally optimal solutions to the Gromov-Wasserstein problem using the sum-of-squares hierarchy. Our main contributions include establishing simplified versions of both Putinar-type and Schm\"udgen-type moment hierarchies that leverage marginal constraints, proving convergence rates of $\mathcal{O}(1/r)$ for these hierarchies, and demonstrating that the resulting SOS-based distance measures satisfy metric properties including the triangle inequality. The theoretical guarantees we provide, particularly the convergence rates and metric properties, establish these hierarchies as principled approximations to the GW problem that are both theoretically sound and computationally tractable. We believe this bridges an important gap between the theoretical understanding of GW distances and practical algorithms for computing them.

Several promising directions remain for future works. A natural extension would be to develop continuous versions of our discrete framework and study the associated gradient flows on the space of metric measure spaces, similar to the theory of distortion distance presented in \cite{Sturm:06}. This could provide new insights into the geometry of metric measure spaces and potentially lead to novel computational methods for the gradient flows. Another important direction is the development of scalable numerical methods for computing these SOS-based distances. While our current approach provides theoretical guarantees, practical implementations for large-scale problems will require careful algorithmic design, possibly incorporating techniques from first-order optimization methods and exploiting the structures of the problem. Such developments would significantly enhance the practical applicability of our framework in machine learning and life-science applications.

\section*{Acknowledgements}

The authors wish to thank Di Hou, Tianyun Tang, and Kim-Chuan Toh for generously sharing their software package RNNAL \cite{hou2025lowrankaugmentedlagrangianmethod}.  The authors also wish to thank the two anonymous reviewers for useful comments that have helped improve the manuscript.

\bibliographystyle{alpha}
\bibliography{gwsos}

\end{document}

%% file: macros.tex
\usepackage{amssymb,amsthm,amsmath}
\usepackage[pdftex]{graphicx}
\usepackage{graphicx}
\usepackage[margin=1in]{geometry}
\usepackage{color}
\usepackage{hyperref}
\usepackage{subcaption}
\usepackage{float}
\usepackage{longtable}
\usepackage{booktabs}
\usepackage{multirow}
\usepackage{makecell}

\newtheorem{theorem}{Theorem}[section]
\newtheorem{lemma}[theorem]{Lemma}
\newtheorem{proposition}[theorem]{Proposition}

\theoremstyle{definition}

\theoremstyle{remark}

\usepackage[textwidth=2.450cm, textsize=tiny]{todonotes}

\newcommand{\CX}{\mathcal{X}}
\newcommand{\CY}{\mathcal{Y}}
\newcommand{\CZ}{\mathcal{Z}}
\newcommand{\muX}{\mu_{\CX}}
\newcommand{\muY}{\mu_{\CY}}

\newcommand{\RR}{\mathbb{R}}

\newcommand{\lb}{\mathrm{lb}}
\newcommand{\mlb}{\mathrm{mlb}}
\newcommand{\NN}{\mathbb{N}}
\newcommand{\M}{\mathbf{M}}

\newcommand{\XX}{\mathbb{X}}
\newcommand{\YY}{\mathbb{Y}}
\newcommand{\ZZ}{\mathbb{Z}}

\newcommand{\iy}{\mathit{z}}
\newcommand{\CT}{\mathcal{T}}
\newcommand{\CQ}{\mathcal{Q}}

\newcommand{\vv}{\mathbf{v}}

\newcommand{\Proj}{\mathbf{P}}
\newcommand{\growbeta}{B}

\newcommand{\muxi}{(\mu_{\CX})_i}
\newcommand{\muyj}{(\mu_{\CY})_j}

\allowdisplaybreaks